\documentclass[12pt]{amsart}

\usepackage{epsfig, color}

\usepackage{amsmath,amssymb}
\usepackage{amsthm}
\usepackage{hyperref}
\usepackage[margin=1.0in]{geometry}

\numberwithin{equation}{section}

\newtheorem{prop}{Proposition}
\newtheorem{lemma}[prop]{Lemma}

\newtheorem{thm}[prop]{Theorem}
\newtheorem{cor}[prop]{Corollary}
\newtheorem{conj}[prop]{Conjecture}
\numberwithin{prop}{section}

\theoremstyle{definition}
\newtheorem{defn}[prop]{Definition}

\newtheorem{rmk}[prop]{Remark}

\newcommand{\del}{\partial}
\newcommand{\delb}{\bar{\partial}}
\newcommand{\dt}{\frac{\partial}{\partial t}}
\newcommand{\brs}[1]{\left| #1 \right|}

\newcommand{\gG}{\Gamma}
\renewcommand{\gg}{\gamma}
\newcommand{\gD}{\Delta}
\newcommand{\gd}{\delta}
\newcommand{\gs}{\sigma}

\newcommand{\gl}{\lambda}

\newcommand{\gw}{\omega}
\newcommand{\ga}{\alpha}
\newcommand{\gb}{\beta}
\renewcommand{\ge}{\epsilon}
\newcommand{\N}{\nabla}
\newcommand{\FF}{\mathcal F}

\newcommand{\KK}{\mathcal K}
\newcommand{\PP}{\mathcal P}

\newcommand{\til}[1]{\widetilde{#1}}

\renewcommand{\bar}[1]{\overline{#1}}

\renewcommand{\i}{\sqrt{-1}}

\newcommand{\hook}{\mathbin{\hbox{\vrule height2.4pt width4.5pt depth-2pt
\vrule height5pt width0.4pt depth-2pt}}}

\newcommand{\bj}{\bar{j}}

\newcommand{\bl}{\bar{l}}

\newcommand{\bgb}{\bar{\beta}}
\newcommand{\bp}{\bar{p}}
\newcommand{\bq}{\bar{q}}

\newcommand{\IP}[1]{\left<#1\right>}

\newcommand{\bga}{\bar{\alpha}}

\DeclareMathOperator{\Sym}{Sym}
\DeclareMathOperator{\Rc}{Rc}
\DeclareMathOperator{\SL}{SL}
\DeclareMathOperator{\Rm}{Rm}

\DeclareMathOperator{\tr}{tr}
\DeclareMathOperator{\Ker}{Ker}

\DeclareMathOperator{\Id}{Id}

\DeclareMathOperator{\kod}{kod}

\DeclareMathOperator{\Nil}{Nil}

\begin{document}

\title[Pluriclosed flow and the geometrization of complex surfaces]{Pluriclosed flow and the geometrization of complex surfaces}

\begin{abstract} We recall fundamental aspects of the pluriclosed flow equation and survey various existence and convergence results, and the various analytic techniques used to establish them.  Building on this, we formulate a precise conjectural description of the long time behavior of the flow on complex surfaces.  This suggests an attendant geometrization conjecture which has implications for the topology of complex surfaces and the classification of generalized K\"ahler structures.
\end{abstract}

\date{\today}

\author{Jeffrey Streets}
\address{Rowland Hall\\
         University of California, Irvine\\
         Irvine, CA 92617}
\email{\href{mailto:jstreets@uci.edu}{jstreets@uci.edu}}

\maketitle

\section{Introduction}

The Enriques-Kodaira classification of compact complex surfaces is a landmark achievement of 20th century mathematics, exploiting deep techniques from complex analysis, partial differential equations and algebraic geometry to give a descriptive classification of complex manifolds of complex dimension $2$.  Despite being dubbed a classification, many questions remain unanswered, and the structure of complex surfaces remains an active area of research to this day.  As the classical uniformization of Riemann surfaces profoundly intertwines complex structures and associated canonical Riemannian geometries, it is natural to try to associate canonical metrics to complex surfaces in order to provide further insight.  The purpose of this article is to describe such a ``geometrization conjecture'' for compact complex surfaces.  Specifically we aim to associate canonical families of Hermitian metrics on all complex surfaces via a universal geometric flow construction, and use properties of the resulting metrics to capture aspects of the underlying complex structure, and in particular yield the classification of Kodaira's Class VII surfaces.  

This strategy for understanding complex surfaces is inspired by Perelman's landmark resolution of Thurston's Geometrization Conjecture for $3$-manifolds using Ricci flow \cite{Perelman1,Perelman2,Perelman3,KL,MorganTian}.  Thurston \cite{Thurston} sought to decompose any compact $3$-manifold into pieces, each of which admits canonical Riemannian metrics, the models of which are inspired by a thorough understanding of possible locally homogeneous spaces.  In 1982, Hamilton introduced the Ricci flow equation and used it to classify compact $3$-manifolds with positive Ricci curvature.  Over the ensuing decades Hamilton further developed the theory of Ricci flow, eventually giving a conjectural description of the singularity formation and long time behavior of the flow, which would lead to a proof of Geometrization.  Exploiting many deep insights into the structure of Ricci flow, Perelman achieved a precise description of this singularity formation, in particular yielding Geometrization.  Amazingly, this one analytic tool provides both the topological decomposition of the underlying $3$-manifold as well as an essentially canonical construction of the relevant geometric structures.

Inspired by this story, the author and Tian sought to bring the philosophy of geometric evolution to bear in understanding the topology and geometry of complex manifolds \cite{HCF, PCF}.  Of course Ricci flow has already had a strong influence in complex, \emph{K\"ahler}, geometry.  The Ricci flow equation preserves K\"ahler geometry, and Cao initiated the study of K\"ahler-Ricci flow, \cite{CaoKRF} using it to reprove the Calabi-Yau \cite{YauCC} and Aubin-Yau Theorems \cite{Aubin, YauCC}, canonically constructing K\"ahler-Einstein metrics on manifolds with $c_1 = 0$, $c_1 < 0$ respectively.  More recently, Song-Tian initiated the analytic minimal model program \cite{SongTian}, seeking to understand the algebraic minimal model program through the singularities of K\"ahler-Ricci flow.  However, there are many examples of complex, \emph{non-K\"ahler} manifolds, starting with the basic example of Hopf who constructed complex structures on $S^3 \times S^1$.  The K\"ahler-Ricci flow cannot be employed on such manifolds, and one can show that for non-K\"ahler metrics the Ricci flow will not even preserve the condition that a metric is Hermitian, thus one must look elsewhere to apply the ideas of geometric evolutions.  Thus the \emph{pluriclosed flow} equation \cite{PCF} was introduced as an extension of K\"ahler-Ricci flow which preserves natural conditions for Hermitian, non-K\"ahler metrics.  

Following the general philosophy of geometric evolutions, one expects the limiting behavior of pluriclosed flow to meaningfully reflect aspects of the complex structure and topology of the underlying surface.  Up to now results have been established which support the main conjecture on the maximal smooth existence time of the flow, and we will recount these below.  Moreover we have phenomenological results in the locally homogeneous setting which suggest the general limiting behavior of the flow in many cases \cite{Boling}.  However, there is not yet a precise conjectural picture of how the flow behaves on most complex surfaces, most importantly Kodaira's Class $\mbox{VII}^+$ surfaces.  These are minimal complex surfaces of negative Kodaira dimension and $b_2 > 0$.  Conjecturally these are all diffeomorphic to $S^3 \times S^1 \# k \bar{\mathbb{CP}}^2$, with a complete description of the different complex structures.  Despite the relative simplicity of the underlying diffeotypes, these represent a large and varied class of surfaces, and their classification remains the main open problem in the Kodaira classification of surfaces.  In this article we formulate precise conjectures on the limiting behavior of the pluriclosed flow on all complex surfaces, and announce a number of new global existence and convergence results in support of these conjectures.

Furthermore, it turns out that the pluriclosed flow is also related to Hitchin's ``generalized geometry'' program.  The first hint of this was shown in \cite{GKRF}, where the author and Tian showed that the pluriclosed flow equation actually preserves the delicate integrability conditions of Gualtieri/Hitchin's ``generalized K\"ahler geometry,'' a subject with roots in mathematical physics.  We will see in part III below that the conjectural framework for pluriclosed flow leads to a classification of generalized K\"ahler structures on complex surfaces.  One example of particular interest is the case of $\mathbb {CP}^2$, where the global existence and convergence of the pluriclosed flow implies a uniqueness theorem for generalized K\"ahler structures, extending Yau's theorem on the uniqueness of the complex structure \cite{YauCCann}.

\vskip 0.2in

\textbf{Acknowledgements:} It is a pleasure to dedicate this article to Professor Gang Tian on the occasion of his 59th birthday.  My collaboration with Tian obviously forms the underpinning of this article, and Tian's many works in K\"ahler geometry helped inspire the picture we paint here.  We thank Georges Dloussky for many helpful conversations.

\vskip 0.2in
\begin{center}
\textbf{\large{Part I: Pluriclosed flow}}
\end{center}
\vskip 0.2in

\section{Existence and basic regularity properties} \label{s:PCFsec}
\subsection{Definition and local existence}

In this section we recount the rudimentary properties of the pluriclosed flow equation.  To begin we define the pluriclosed condition.

\begin{defn} Let $(M^{2n}, g, J)$ be a complex manifold with Hermitian metric $g$, and associated K\"ahler form $\gw(X,Y) = g(X, J Y)$.  The metric $g$ is \emph{K\"ahler} if
\begin{align*}
d \gw = 0.
\end{align*}
The metric $g$ is \emph{pluriclosed} if
\begin{align*}
\i \del \delb \gw = d d^c \gw = 0,
\end{align*}
where $d^c \omega = \i (\delb - \del) \omega = - d \gw(J, J, J)$.  As they are equivalent notions we will often refer to the K\"ahler form $\gw$ as being either K\"ahler or pluriclosed.
\end{defn}

The K\"ahler condition is the simplest, and strongest, integrability condition for a Hermitian metric.  The pluriclosed condition is essentially the only weakening of the K\"ahler condition which is linear in the K\"ahler form.  As the aim is to understand all complex surfaces through the geometry of pluriclosed metrics, the following result of Gauduchon is of fundamental importance.

\begin{thm} \label{t:gauduchon} (\cite{Gauduchon1form}) Let $(M^{2n}, g, J)$ be a compact complex manifold with Hermitian metric $g$.  There exists a unique $u \in C^{\infty}(M)$ such that
\begin{align*}
\i \del \delb \left( e^{2 u} \gw \right)^{n-1} = 0, \qquad \int_M u dV_g = 0.
\end{align*}
\end{thm}
\noindent In particular, applying this theorem in the case $n=2$, we see that every compact complex surface admits pluriclosed metrics.

It is well known that for a K\"ahler manifold, the Levi-Civita connection preserves $J$, and is the unique Hermitian connection associated to the pair $(g,J)$.  When the metric is not K\"ahler, there is a natural one-parameter family of Hermitian connections associated to $(g,J)$ (\cite{Gauduchonconn}).  Two of these are particularly relevant to pluriclosed flow, namely the Bismut and Chern connections, defined by:
\begin{align*}
\IP{\N^B_X Y, Z} =&\ \IP{\N^{LC}_X Y, Z} + \tfrac{1}{2} d^c \omega(X,Y,Z),\\
\IP{\N^C_X Y, Z} =&\ \IP{\N^{LC}_X Y, Z} + \tfrac{1}{2} d \omega(J X,Y,Z).
\end{align*}
These connections induce curvature tensors $\Omega^B$ and $\Omega^C$, and also connections on the canonical bundle, yielding in turn representatives of the first Chern class:
\begin{align*}
\rho^B_{\ga \gb} =&\ g^{\bj i} \Omega^B_{\ga \gb i \bj}, \qquad \rho^C_{\ga \gb} = g^{\bj i} \Omega^C_{\ga \gb i \bj}.
\end{align*}
By Chern-Weil theory we know that $d \rho^{B,C} = 0$, and $\rho^{B,C} \in c_1$.  However, it is not the case in general that $\rho^B$ is of type $(1,1)$.  Of course in the K\"ahler setting $\rho^B = \rho^C = \rho$, the Ricci form of the underlying metric.

The classical quest for canonical metrics on manifolds often centers around existence questions for Einstein metrics, i.e. solutions of
\begin{align*}
\Rc = \gl g
\end{align*}
for some constant $\gl$.  Hamilton \cite{Hamilton3folds} introduced the Ricci flow
\begin{align*}
\dt g =&\ -2 \Rc
\end{align*}
as a tool for constructing such metrics by the parabolic flow method.  Cao \cite{CaoKRF} observed that the Ricci flow equation will preserve the K\"ahler condition, and can be expressed in terms of the K\"ahler form as
\begin{align} \label{f:KRF}
\dt \gw =&\ - \rho = \i \del \delb \log \det g,
\end{align}
where the last expression holds in local complex coordinates.

Given the vast array of analytic tools and structure available for Ricci flow, one would want to use it to understand complex manifolds beyond the K\"ahler setting.  However, for a general Hermitian metric the Ricci tensor is not of $(1,1)$ type and thus the Ricci flow equation will not preserve the class of Hermitian metrics.  For this reason we are forced to define a new equation if we are to preserve aspects of Hermitian, non-K\"ahler geometry.  In \cite{HCF} the author and Tian introduced a family of parabolic equations for Hermitian metrics on complex manifolds.  Let $S_{i \bj} = g^{\bl k} \Omega^C_{k \bl i \bj}$, which is a kind of Ricci tensor defined using the Chern connection, and is always a $(1,1)$ form.  Furthermore, we let $T$ denote the torsion of the Chern connection, and let $Q = T \star T$ denote an \emph{arbitrary} quadratic in $T$ which satisfies $Q \in \Lambda^{1,1}$.  In \cite{HCF} the author and Tian defined the class of flow equations:
\begin{gather}
\dt \gw = - S + Q,
\end{gather}
referring to these generally as \emph{Hermitian curvature flow}.  We showed that for an arbitrary choice of $Q$ this equation is strictly parabolic and satisfies various natural analytic conditions such as smoothing estimates and stability near K\"ahler-Einstein metrics.  In \cite{HCF} a particular choice of $Q$ was identified corresponding to the Euler equation of a certain Hilbert-type functional in Hermitian geometry.  More recently Ustinovskiy \cite{Yurypreprint} showed that for a different choice of $Q$ the flow will preserve various curvature positivity conditions, leading to extensions of the classical Frankel conjecture into non-K\"ahler geomery.  Given the rich diversity of Hermitian geometry, it is natural to expect that different flows, i.e. different choices of $Q$, would be needed to address different situations.  In \cite{PCF} the author and Tian identified a particular choice of $Q$ which yields a flow which preserves the pluriclosed condition, specifically
\begin{align*}
Q^1_{i \bj} = g^{\bl k} g^{\bq p} T_{i k \bp} T_{\bj \bl q}.
\end{align*}
As it turns out the resulting flow equation has several different useful manifestations, and we record the fundamental definition including some of these forms.

\begin{defn} Let $(M^{2n}, J)$ be a complex manifold.  We say that a one-parameter family of pluriclosed metrics $\gw_t$ is a solution of \emph{pluriclosed flow} if
\begin{align*}
\dt \gw =&\ - S + Q^1.
\end{align*}
This is an example of Hermitian curvature flow, which is well posed for arbitrary Hermitian metrics.  With pluriclosed initial data it can also be expressed using the Bismut curvature as
\begin{align*}
\dt \gw =&\ - \rho_B^{1,1}.
\end{align*}
It is useful to furthermore express the flow in local complex coordinates, yielding
\begin{align} \label{f:PCF}
\dt \gw =&\ \del \del^*_{\gw} \gw + \delb \delb^*_{\gw} \gw + \i \del \delb \log \det g.
\end{align}
\end{defn}

As an example of Hermitian curvature flow, there exist short time solutions on compact manifolds as remarked above.  Arguing specifically in this setting, it turns out that the operator $\rho_B^{1,1}$, restricted to the class of pluriclosed metrics, is strictly elliptic.  In fact, the symbol of the linearized operator is just the Laplacian with respect to the given metric.  This renders pluriclosed flow a strictly parabolic equation, and so by appealing to general theory we can obtain short time existence on compact manifolds.  Moreover, comparing (\ref{f:PCF}) and (\ref{f:KRF}), it is natural to expect that if we start the flow with K\"ahler initial data then it is a solution of K\"ahler-Ricci flow.  This in fact holds, so to summarize:

\begin{thm} (cf. \cite{PCF} Theorem 1.2) \label{t:ste} Let $(M^{2n}, J)$ be a compact complex manifold.  Suppose $\gw_0$ is a pluriclosed metric on $M$.  Then there exists $\ge > 0$ and a unique solution to pluriclosed flow $\gw_t$ with initial condition $\gw_0$.  If $\gw_0$ is K\"ahler, then $\gw_t$ is the unique solution to K\"ahler-Ricci flow with initial data $\gw_0$.
\end{thm}

Having obtained short time solutions, the main task is to describe the maximal smooth existence time, as well as the limiting behavior at this time, in general a formidable task.  The remainder of this paper is devoted to giving precise conjectures on these questions in the specific case of complex surfaces.

\begin{rmk} We point out that a different geometric flow approach has been introduced to study complex surfaces, the Chern-Ricci flow:
\begin{align*}
\dt \gw =&\ - \rho^C.
\end{align*}
This flow also preserves the pluriclosed condition, and reduces to a scalar PDE modeled on the parabolic complex Monge-Ampere equation.  A crucial qualitative distinction between this flow and pluriclosed flow is that for Chern-Ricci flow the torsion is fixed as a tensor along the flow, i.e. $d \gw_t = d \gw_0$, whereas for pluriclosed flow the torsion tensor satisfies a parabolic PDE.  

Despite this difference there are similarities between the results and expectations of pluriclosed flow and Chern-Ricci flow, especially in the cases of positive Kodaira dimension, where it is natural to rescale by ``blowing down'' the flow, which then yields $d \gw_t = e^{-t} d \gw_0$, allowing for convergence to K\"ahler metrics.  Specifically, the analogue of Conjectures \ref{c:gentypeconj} and \ref{c:propellipconj} below were shown to hold for Chern-Ricci flow (\cite{WhineCRF1} Theorem 1.7, \cite{WhineCRF2} Theorem 1.1).  

Differences start to appear in Kodaira dimension zero.  For instance, as shown in \cite{GillCRF}, given an arbitrary Hermitian metric on the torus, the Chern-Ricci flow will exist globally and converge to a \emph{Chern-Ricci flat}, but not necessarily \emph{flat}, metric.  Here the fact that $d \gw_t = d \gw_0$ prevents the flow converging to a K\"ahler metric when starting from a non-K\"ahler metric in this setting.  This is related to the fact that there is an infinite dimensional moduli space of pluriclosed Chern-Ricci flat metrics on the torus, using perturbations of the flat metric via $\delb \ga + \del \bga$.  Alternatively, Theorem \ref{t:torustheorem} shows that the pluriclosed flow on the torus, with arbitrary initial data, exists globally and converges to a flat metric.

The differences become even more stark for Kodaira dimension $-\infty$, where for instance pluriclosed flow has fixed points on Hopf surfaces, whereas Chern-Ricci flow always encounters a finite time singularity.  Also on $\mathbb{CP}^2$ the Fubini-Study metric is stable for normalized pluriclosed flow (also see Theorem \ref{t:CP2} below for more general convergence results), whereas normalized Chern-Ricci flow will satisfy $d \gw_t = e^t d \gw_0$, and so cannot converge to Fubini-Study for non-K\"ahler initial data.  In a strange quirk, the geometry of Inoue surfaces again requires a blowdown to obtain a geometric limit, so here again the expectations between the two flows agree, and the natural analogue of Conjecture \ref{c:Inoueconj} for Chern-Ricci flow was shown to hold for certain initial data \cite{CRFonInoue}.
\end{rmk}

\subsection{Pluriclosed flow as a gradient flow} \label{ss:gradient}

The original motivation for defining pluriclosed flow came from complex geometry, aiming to preserve natural properties of Hermitian metrics.  As it turns out, this equation has a remarkable and useful relationship to the Ricci flow or, more precisely, the generalized Ricci flow which couples to the heat equation for a closed three-form.

\begin{thm} \label{t:PCFRGF} (cf. \cite{PCFReg} Theorem 6.5) Let $(M^{2n}, \gw_t, J)$ be a solution to pluriclosed flow.  Let $(g_t, H_t = d^c \gw_t)$ be the associated $1$-parameter families of Riemannian metrics and Bismut torsion forms.  Then
\begin{gather}
\begin{split}
\dt g =&\ - 2 \Rc + \tfrac{1}{2} H^2 - L_{\theta^{\sharp}} g,\\
\dt H =&\ \gD_d H - L_{\theta^{\sharp}} H.
\end{split}
\end{gather}
\end{thm}

Notice that we may apply a family of diffeomorphisms to remove the Lie derivative terms, yielding a solution of the system of equations
\begin{gather} \label{f:GRF}
\begin{split}
\dt g =&\ - 2 \Rc + \tfrac{1}{2} H^2,\\
\dt H =&\ \gD_d H.
\end{split}
\end{gather}
This system of equations originally arose in mathematical physics in the context of renormalization group flow in the theory of $\gs$-models.  The author studied this system under the name ``connection Ricci flow'' \cite{Streetsexpent} and ``generalized Ricci flow'' \cite{StreetsTdual} due to the relationship to the curvature of the Bismut connection and generalized geometry.  Notice that in the context of pluriclosed flow, if we apply the relevant gauge transformation to produce a solution $(g_t, H_t)$ to (\ref{f:GRF}), then $g_t$ remains a pluriclosed metric, but with respect to the appropriately gauge-modified family of complex structures.  This seemingly minor point is actually an essential and rich feature of the generalized K\"ahler-Ricci flow, explained in \S \ref{s:GKRF}.  But first, the primary analytic consequence of Theorem \ref{t:PCFRGF} is the realization of pluriclosed flow as a gradient flow.  As shown in \cite{OSW}, (\ref{f:GRF}) is the gradient flow of the first eigenvalue of a certain Schr\"odinger operator, extending Perelman's monotonicity formulas for Ricci flow.  We briefly describe this construction, and the reader should consult (\cite{OSW,Perelman1}) for further detail.  

To begin we define the functional
\begin{align} \label{f:Fdef}
\FF(g,H,f) = \int_M \left( R - \tfrac{1}{12} \brs{H}^2 + \brs{\N f}^2 \right) e^{-f} dV_g,
\end{align}
for $f \in C^{\infty}(M)$.  This functional obeys a monotonicity formula when $f$ obeys the appropriate conjugate heat equation, which is the adjoint of the heat equation with respect to the spacetime $L^2$ norm, in particular
\begin{align} \label{f:BHE}
\dt f =&\ - \gD f - R + \tfrac{1}{4} \brs{H}^2.
\end{align}
Given $(g_t, H_t)$ a solution of generalized Ricci flow, and $f_t$ is an associated solution of (\ref{f:BHE}), then one obtains the equation
\begin{align} \label{f:Fmon}
\frac{d}{dt} \FF(g_t,H_t,f_t) =&\ \int_M \left[ 2 \brs{\Rc - \tfrac{1}{4} H^2 + \N^2 f}^2 + \brs{ d^* H + \N f \hook H}^2\right] e^{-f} dV_g.
\end{align}
The right hand side is of course nonnegative, so that $\FF$ is monotone nondecreasing.  Underlying this monotonicity is the fact that generalized Ricci flow is the gradient flow of the lowest eigenvalue $\gl$ of the operator $-4 \gD + R - \tfrac{1}{12} \brs{H}^2$, characterized via
\begin{align} \label{f:ldef}
\gl(g,H) := \inf_{\left\{ f |\ \int_M e^{-f} dV_g = 1 \right\}} \FF(g,H,f).
\end{align}

\begin{thm} \label{t:OSW} (\cite{OSW} Proposition 3.4) Generalized Ricci flow is the gradient flow of $\gl$.
\end{thm}

These monotonicity formulas have strong consequences for the singularity formation and conjectural framework for pluriclosed flow.  In particular, this monotonicity formula motivates a key concept, that of a \emph{generalized Ricci soliton}.  These are triples $(g,H,f)$ such that
\begin{gather}
\begin{split}
\Rc - \tfrac{1}{4} H^2 + \N^2 f =&\ 0\\
d^* H + \N f \hook H =&\ 0,
\end{split}
\end{gather}
which by (\ref{f:Fmon}) correspond to critical points for $\FF$ or $\gl$.  Given such a triple $(g,H,f)$, the solution to generalized Ricci flow with this initial data evolves via pullback by the $1$-parameter family of diffeomorphisms generated by $- \N f$, and so these solutions are self-similar, and thus a natural, more general notion of a fixed point for generalized Ricci flow.

\section{Conjectural existence properties}

In this section we describe the conjectural maximal smooth existence time for pluriclosed flow on compact complex manifolds, and give a more refined picture in the case of complex surfaces, following \cite{PCFReg}.  We first recall the fundamental theorem of Tian-Zhang \cite{TianZhang} on the maximal smooth existence of K\"ahler-Ricci flow.  We then extend the formal definition behind this result to the case of pluriclosed flow, and state the relevant conjectures.  To finish we recall a result of \cite{PCFReg} giving a characterization of the relevant positive cone in cohomology on complex surfaces which allows us to explicitly compute the formal existence time for the flow in part II.

\subsection{Sharp local existence for K\"ahler-Ricci flow}

To understand the formal picture of the long time existence 
and singularity formation of K\"ahler-Ricci flow, we first study the flow at the 
level of cohomology.

\begin{defn} \label{d:Dolbeault} Let $(M^{2n}, J)$ be a K\"ahler manifold.  Let
\begin{align*}
H^{1,1}_{\mathbb R} := \frac{ \left\{ \Ker d : \Lambda^{1,1}_{\mathbb R} \to \Lambda^{3} 
\right\}}{ \left\{\i \del \delb f\ |\ f \in C^{\infty} \right\} }.
\end{align*}
Furthermore, define the 
\emph{K\"ahler cone} via
\begin{align*}
\KK := \left\{ [\psi] \in H^{1,1}\ |\ \exists\ \gw \in [\psi],\ \gw > 0 
\right\}.
\end{align*}
\end{defn}

The structure of the K\"ahler cone plays a fundamental role in understanding 
the singularity formation of the K\"ahler-Ricci flow.  First note that an elementary consequence of the K\"ahler-Ricci flow equation is that
\begin{align} \label{f:Kflowclassevol}
[\gw_t] = [\gw_0] - t c_1.
\end{align}
Thus the K\"ahler class moves along a ray in $\KK$, and we thus obtain an upper bound for the possible smooth existence time:

\begin{lemma} \label{l:KflowtimeUB} Let $(M^{2n}, g_0, J)$ be a
K\"ahler manifold.  Let
\begin{align} \label{f:Ktaustardef}
\tau^*(\gw_0) := \sup \{ t \geq 0\ |\ [\gw_0] - t c_1 \in \KK \},
\end{align}
and let $T$ denote the maximal smooth existence time for the K\"ahler-Ricci flow
with initial condition $g_0$.  Then $T \leq \tau^*(\gw_0)$.
\begin{proof} Let $\gw_t$ denote the one parameter family of K\"ahler forms 
evolving by K\"ahler-Ricci flow with initial condition $\gw_0$.  If the flow 
existed smoothly for some time $t > \tau^*$, then in particular by (\ref{f:Kflowclassevol}) there exists a 
smooth positive definite metric in $[\gw_0] - t c_1$, contradicting the 
definition of $\tau^*$.
\end{proof}
\end{lemma}

Informally this lemma states that the flow must go singular by the time the associated family of K\"ahler classes leaves the K\"ahler cone.  In view of this, the natural question to ask is if singularities can possibly form without leaving the K\"ahler cone.  The answer is no, due to Tian-Zhang, meaning that $\tau^*(\gw_0)$ is the maximal smooth existence time of the flow.

\begin{thm} \label{t:TianZhang} (\cite{TianZhang} Proposition 1.1) Let $(M^{2n}, g_0, J)$ be a compact K\"ahler manifold.  The maximal smooth solution of K\"ahler-Ricci flow with initial condition $g_0$ exists on $[0, \tau^*(\gw_0))$.
\end{thm}

The proof requires the development of a priori estimates for the metric along the flow.  The fundamental role played by the formal considerations above is that, for times $t < \tau^*(\gw_0)$, it is possible to reduce the K\"ahler-Ricci flow to a scalar PDE modeled on the parabolic complex Monge-Ampere equation.  This allows for various delicate applications of the maximum principle to obtain control over the metric as long as $t < \tau^*(\gw_0)$.

\subsection{A positive cone and conjectural existence for pluriclosed flow}

Now we follow the discussion of the previous subsection and investigate the formal existence time for pluriclosed flow.  First we define the relevant positive cone, this time in an Aeppli cohomology space.

\begin{defn} Let $(M^{2n}, J)$ be a complex manifold.  Define the \emph{real $(1,1)$ Bott-Chern cohomology} via
\begin{align*}
H^{1,1}_{BC, \mathbb R} := \frac{\left\{ \Ker d : 
\Lambda^{1,1}_{\mathbb R} \to \Lambda^{2,2}_{\mathbb R}\right\}}{\left\{\i \del \delb f |\ f \in C^{\infty} \right\}}.
\end{align*}
Furthermore, define the 
\emph{real $(1,1)$ Aeppli cohomology} via
\begin{align*}
H^{1,1}_{\del + \delb, \mathbb R} := \frac{\left\{ \Ker \i \del\delb : 
\Lambda^{1,1}_{\mathbb R} \to \Lambda^{2,2}_{\mathbb R}\right\}}{\left\{\del 
\bga + \delb \ga\ |\ \ga \in \Lambda^{1,0} \right\}}.
\end{align*}
The restriction to real $(1,1)$ forms make these different spaces from what is usually referred to as Bott-Chern and Aeppli cohomology.  Nonetheless in what follows we will drop the $\mathbb R$ from the notation and always mean these spaces defined above.  Lastly, we define the 
\emph{$(1,1)$ Aeppli positive cone} via
\begin{align*}
\PP := \left\{ [\psi] \in H^{1,1}_{\del + \delb}\ |\ 
\exists\ \gw \in [\psi],\ \gw > 0 \right\}.
\end{align*}
\end{defn}

Note that $\PP$ consists precisely of the $(1,1)$ Aeppli classes represented by pluriclosed metrics.  Similarly to (\ref{f:Kflowclassevol}), we want to derive an ODE for the Aeppli classes associated to a solution of pluriclosed flow.  First note that, on a general complex manifold, the first Chern class $c_1$ is an element of $(1,1)$ Bott-Chern cohomology.  There is a natural inclusion map $i : H^{1,1}_{BC} \to H^{1,1}_{\del + \delb}$, and using this inclusion we consider $c_1$ as an element of $(1,1)$ Aeppli cohomology, without further notation.  Thus, given a solution $\gw_t$ of pluriclosed flow (\ref{f:PCF}), we compute as an equation of $(1,1)$ Aeppli classes,
\begin{align*}
\frac{d}{dt} [\gw_t] = - [ \rho_B^{1,1}] = [\del \del^*_{\gw_t} \gw_t + \delb \delb^*_{\gw_t} \gw_t - \rho_C(\gw_t)] = - c_1.
\end{align*}
Thus precisely as in the K\"ahler-Ricci flow case we obtain, as an equation of $(1,1)$ Aeppli classes,
\begin{align} \label{f:PCFODE}
[\gw_t] = [\gw_0] - t c_1.
\end{align}
This formal calculation allows us to derive an upper bound for the maximal smooth existence time of a solution to pluriclosed flow, whose proof is identical to that of Lemma \ref{l:KflowtimeUB}.

\begin{lemma} \label{l:flowtimeUB} Let $(M^{2n}, g_0, J)$ be a complex manifold 
with pluriclosed metric.  Let
\begin{align} \label{f:taustardef}
\tau^*(\gw_0) := \sup \{ t \geq 0\ |\ [\gw_0] - t c_1 \in \PP \},
\end{align}
and let $T$ denote the maximal smooth existence time for the pluriclosed flow 
with initial condition $g_0$.  Then $T \leq \tau^*(\gw_0)$.
\end{lemma}

What follows is the main conjecture guiding the study of pluriclosed flow.  While Lemma \ref{l:flowtimeUB} indicates the elementary fact that the maximal existence time for the flow can be no larger than $\tau^*(\gw_0)$, Conjecture \ref{c:mainflowconj} indicates that the flow is actually smooth up to time $\tau^*(\gw_0)$, i.e. that it equals the maximal existence time.
\begin{conj}[\textbf{Main Existence Conjecture}]\label{c:mainflowconj} Let $(M^{2n}, g_0, J)$ be a compact complex manifold with pluriclosed metric.  The maximal smooth solution of pluriclosed flow with initial condition $g_0$ exists on $[0, \tau^*(\gw_0))$.
\end{conj}
\noindent This conjecture first appeared in our joint work with Tian (\cite{PCFReg} Conjecture 5.2), inspired by Theorem \ref{t:TianZhang}.

\subsection{Characterizations of positive cones}

As a guide for the nature of singularity formation as one leaves the positive 
cone it is useful to have a characterization of the necessary and sufficient 
conditions for cohomology classes to lie in the appropriate positive cone.  We give such a characterization in the case $n=2$ in this subsection.  First, for a given complex surface define
\begin{align*}
\Gamma = \frac{\{d a \in \Lambda^{1,1}_{\mathbb R} \}}{\left\{ \i \del \delb f| f \in C^{\infty} \right\}}.
\end{align*}
By (\cite{Telemancone} Lemma 2.3), This is identified with a subspace of $\mathbb R$, via the $L^2$ inner product 
with a pluriclosed metric $\omega$.  Note that if $(M^4, J)$ is 
K\"ahler then the $\del\delb$-lemma holds and so $\gG = \{0\}$.   Further arguments of
(\cite{Telemancone} Lemma 2.3) in fact show that the vanishing of $\gG$ implies 
the manifold is K\"ahler.  Thus, on a non-K\"ahler surface we may choose a 
positive generator $\gg_0$ for $\gG$, and since the space of
pluriclosed metrics on $M$ is connected, this orientation is well-defined.  
This form $\gg_0$ plays a key role in the characterization of the positive cone 
$\PP$ in the next theorem.

\begin{thm} \label{t:conecharacterization} (\cite{PCFReg} Theorem 5.6) Let $(M^4, J)$ be a
complex non-K\"ahler surface.
 Suppose $\phi \in
\Lambda^{1,1}$ is pluriclosed.  Then $[\phi] \in \mathcal P$ if
and
only if
\begin{enumerate}
\item{$\int_M \phi \wedge \gamma_0 > 0$}
\item{$\int_D \phi > 0 \mbox{ for every effective divisor with negative self
intersection}$.}
\end{enumerate}
\end{thm}

A natural question is whether there is a characterization 
of $\PP$ in higher dimensions.  In such cases it is not even 
clear how to define natural conditions which only depend on Aeppli cohomology 
classes.  As the quantity $\int_M \gw \wedge \gg_0$ will be fixed along a solution to pluriclosed flow, as a corollary of this theorem we obtain a clean characterization of the quantity $\tau^*$:

\begin{cor} \label{c:conecor} Let $(M^4, J)$ be a complex non-K\"ahler surface.  Given $\gw_0$ a pluriclosed metric, one has
\begin{align*}
\tau^*(\gw_0) = \sup \left\{t \geq 0\ |\ \int_D \gw_0 - t c_1 > 0 \mbox{ for } D^2 < 0 \right\}.
\end{align*}
\end{cor}

\section{$(1,0)$-form reduction}

As described above, the proof of Theorem \ref{t:TianZhang} rests on the key fact that the K\"ahler-Ricci flow, for times $t < \tau^*(\gw_0)$, the K\"ahler-Ricci flow can be reduced to a scalar equation.  Due to the $\del\delb$-lemma, we know that locally K\"ahler metrics admit scalar potential functions, and thus one expects any variation of K\"ahler metrics to reduce to a variation of K\"ahler potentials, up to global topological considerations.  Locally this PDE is modeled on the parabolic complex Monge-Ampere equation,
\begin{align} \label{f:LMA}
\frac{\del f}{\del t} =&\ \log \det \i \del \delb f,
\end{align}
where $\gw = \i \del \delb f$.  Using the scalar reduction, various maximum principle arguments are employed to obtain $C^{\infty}$ estimates for $f$, and thus the K\"ahler-Ricci flow (\cite{TianZhang}).

Turning to the pluriclosed flow, we first note that in non-K\"ahler geometry, Hermitian metrics, even pluriclosed, cannot be described by a single potential function.  Thus it is not reasonable to expect to reduce the pluriclosed flow to a scalar PDE.  Instead, pluriclosed metrics admit local potential $(1,0)$ forms.  In particular, a short agument using the Dolbeault lemma shows that for any pluriclosed metric, locally there exists a $(1,0)$-form $\ga$ such that $\gw = \delb \ga + \del \bga$.  Thus for a solution to pluriclosed flow, locally we can express $\gw_t = \gw_{\ga_t} := \delb \ga_t + \del \bga_t$, and using (\ref{f:PCF}) one can show that $\ga_t$ should locally satisfy the PDE
\begin{align} \label{f:1formred}
\frac{\del \ga}{\del t} =&\ \delb^*_{g_{\ga}} \gw_{\ga} - \frac{\i}{2} \del \log \det g_{\ga}.
\end{align}
Observe that this is a strict generalization of (\ref{f:LMA}), where if the metric is K\"ahler and $\ga_t = \frac{\i}{2} \del f$, then the PDE for $\ga$ corresponds to that satisfied by the gradient of a function evolving by (\ref{f:LMA}).  As it turns out, equation (\ref{f:1formred}) is degenerate parabolic, with the degeneracy arising from the redundancy wherein $\ga$ and $\ga + \del f$ describe the same metric for $f \in C^{\infty}(M,\mathbb R)$.

A full, positive resolution of Conjecture \ref{c:mainflowconj} will hinge on a complete understanding of equation (\ref{f:1formred}).  Using this reduced equation, we have acheived many global existence and convergence results which confirm Conjecture \ref{c:mainflowconj} in a variety of cases, and these are described below.  We will not give a full account of the proofs of these results here, but rather describe one key estimate underlying all of these proofs.  In establishing regularity of K\"ahler-Ricci flow, it is crucial to obtain the $C^{2,\ga}$ estimate for the potential in the presence of $C^{1,1}$ estimates.  This can be achieved using Calabi/Yau's $C^3$ estimate (\cite{Calabiaffine} \cite{YauCC}) or the Evans-Krylov style of estimates (\cite{EvansC2a, Krylov}).  On the other hand, thinking in terms of the metric tensor, this is a $C^{\ga}$ estimate in the presence of uniform parabolicity bounds.  As the pluriclosed flow is a parabolic system of
equations for the Hermitian metric $g$, an estimate of this kind would be analogous to the
DeGiorgi-Nash-Moser/Krylov-Safonov \cite{DeGi,Nash,Moser,KS1,KS2} estimate for
uniformly parabolic equations.  However, these results are
false for general \emph{systems} of equations \cite{DeGiorgiCE}.  Despite these challenges we are able to obtain this an estimate of this kind by uncovering a kind of convexity structure described below.

We give an informal statement of the result here, which combines (\cite{StreetsPCFBI} Theorems 1.7, 1.8), referring there for the precise claims.

\begin{thm} \label{EKthm} Let $(M^{2n}, J)$ be a compact complex manifold.
 Suppose $g_t$ is a solution to the pluriclosed flow on $[0,\tau)$.  Suppose there
exist a constant $\gl > 0$ such that
\begin{align*}
\gl g_0 \leq g_t.
\end{align*}
Then there exist uniform $C^{\infty}$ estimates for $g_t$ on $[0,\tau)$.
\end{thm}

The crucial point behind this theorem is a sharp differential inequality for a delicate combination of first derivatives of the potential $\ga$.  In particular, consider the section of $\Sym^2(T \oplus T^*)$ defined by
\begin{align*}
W = \left(
\begin{matrix}
g - \del \ga g^{-1} \delb \bga & \i \del \ga g^{-1}\\
\i g^{-1} \delb \bga & g^{-1}
\end{matrix}
\right)
\end{align*}
The form of this matrix is inspired in part by Legendre transforms \cite{StreetsWarren}, and in part from its appearance as the ``generalized metric,'' in generalized geometry \cite{GualtieriThesis}, where $i \del \ga$ plays the role of the ``$B$-field.''
A delicate computation shows that $L W \leq 0$, where $L$ is the Laplacian of the time-dependent metric.  Using this together with the fact that $\det W = 1$ it is possible to obtain a $C^{\ga}$ estimate for the matrix $W$, which yields $C^{\ga}$ estimates for $g$.  Via a blowup argument one then obtains all higher order estimates.

Thus using this theorem we see that the remaining obstacle to establish Conjecture \ref{c:mainflowconj} is to obtain a uniform lower bound for the metric along the flow.  This is a significant hurdle, which is overcome in the K\"ahler setting through a delicate combination of maximum principles for the parabolic complex Monge-Ampere equation.  There are many cases where Conjecture \ref{c:mainflowconj} can be established by exploiting some further underlying geometric structures.  The simplest case to address is that of manifolds with nonpositive bisectional curvature:

\begin{thm} \label{t:LTE} (\cite{StreetsPCFBI} Theorem 1.1) Suppose $(M^4, J, h)$ is a compact complex surface, with Hermitian metric $h$ with nonpositive holomorphic bisectional curvature.  Given $g_0$ a pluriclosed metric on $M$, the solution to normalized pluriclosed flow exists for all time.
\end{thm}

In general the limiting behavior at infinity in the above theorem can be delicate, involving contraction of divisors as well as collapsing.  With a further curvature restriction we can obtain a convergence result as well:

\begin{thm} \label{t:negcurvthm} (\cite{StreetsPCFBI} Theorem 1.1) Suppose $(M^4, J, h)$ is a compact complex surface, with Hermitian metric $h$ with constant negative bisectional curvature.  Given $g_0$ a pluriclosed metric on $M$, the solution to normalized pluriclosed flow exists for all time and converges to $h$.
\end{thm}

Also, we can obtain convergence to a flat metric in the case of tori, which requires further a priori estimates and use of the $\FF$-functional described in \S \ref{ss:gradient}.

\begin{thm} \label{t:torustheorem} (\cite{StreetsPCFBI} Theorem 1.1) Let $(M^4, J)$ be biholomorphic to a torus.  Given $g_0$ a pluriclosed metric, the solution to pluriclosed flow with initial data $g_0$ exists on $[0,\infty)$ and converges to a flat K\"ahler metric.
\end{thm}

This theorem confirms the basic principle that pluriclosed flow cannot develop ``local'' singularities.  For instance, on any manifold one can choose initial metrics whereby the Ricci flow develops topologically trivial neckpinches.  The rigidity of pluriclosed metrics/pluriclosed flow apparently prevents the construction of such singular solutions.

\section{Pluriclosed flow of locally homogeneous surfaces}

Having described fundamental analytic aspects governing the regularity of pluriclosed flow, we now turn to describing the behavior of locally homogeneous metrics, which in principle give the prototypical behavior for different classes of complex surfaces.  As we will see the pluriclosed flow naturally incorporates two classical points of view on canonical geometries: ``geometric structures'' in the sense of Thurston and (K\"ahler) Einstein metrics.  In this section we recall fundamental aspects of Thurston geometries and how pluriclosed flow of locally homogeneous metrics relates to these structures.

\subsection{Wall's Classification}

In 1985 Wall, directly inspired by the Thurston geometrization conjecture \cite{Thurston}, sought to understand complex surfaces through the use of Thurston's ``model geometries.''  By combining elements of Lie theory with results from the Kodaira classification he gave a complete classification of complex surfaces admitting model geometric structure.  Following \cite{Wall} we will briefly recall these results and later describe the relationship of pluriclosed flow to this classification.  We recall the fundamental definition:

\begin{defn} \label{d:models} A \emph{model geometry} is a triple $(X, g, G)$ such that $(X, g)$ is a complete, simply connected Riemannian manifold, $G$ is a group of isometries acting transitively on $(X, g)$, and $G$ has a discrete subgroup $\gG$ such that $\gG \backslash X$ is compact.
\end{defn}

Let us very briefly describe the model geometries up to dimension $4$, only indicating the space $X$.  In dimension $1$ there is a unique geometry, the real line.  In dimension $2$ there are three, namely $S^2, \mathbb R^2,$ and $\mathbb H^2$.  In dimension three there are eight geometries, classified by Thurston \cite{ThurstonBook}.  First there are $S^3, \mathbb R^3$, and $\mathbb H^3$, as well as the product geometries $\mathbb S^2 \times \mathbb R$ and $\mathbb H^2 \times \mathbb R$.  Another example is $\til{SL}_2$, the universal cover of the unit tangent bundle of $\mathbb H^2$, with metric left-invariant under the natural Lie group structure.  Also one has $Nil^3$, the unique simply connected nilpotent Lie group in three dimensions, as well as a solvable Lie group $Sol^3$, realized as an extension $\mathbb R^2 \ltimes_{\ga} \mathbb R$, where
$\ga(t)(x,y) = (e^t x, e^{-t} y)$, again with left invariant metric.

In four dimensions these geometries were classified by Filipkiewicz \cite{FilipThesis}. First there are the irreducible symmetric spaces $S^4, \mathbb H^4, \mathbb {CP}^2, \mathbb {CH}^2$, as well as products of all lower dimensional examples above.  There is a four dimensional nilpotent Lie group $Nil^4$ and a family of solvable Lie groups $Sol^4_{m,n}$ which we will not describe as these do not admit compatible complex structures.  Another class of solvable Lie groups arises, $Sol_0^4 = \mathbb R^3 \ltimes_{\gd} \mathbb R$, where $\gd(t)(x,y,z) = (e^t x, e^t y, e^{-2t} z)$.  Lastly one has $Sol_1^4$ which is the Lie group of matrices
\begin{align*}
\left\{ \left(
\begin{matrix}
1 & b & c\\
0 & \ga & a\\
0 & 0 & 1
\end{matrix}\right)\ :\ \ga,a,b,c \in \mathbb R, \ga > 0 \right\}.
\end{align*}
Another geometry denoted $F^4$ is identified which we ignore as it admits no compact models.

In describing the existence of compatible complex structures on these geometries, some subtleties arise.  We recall the main theorems of Wall \cite{Wall} here.

\begin{thm} \label{t:Wallthm}
\begin{enumerate}
\item A model geometry $X$ carries a complex structure compatible with the maximal connected group of isometries if and only if $X$ is one of:
\begin{align*}
\mathbb {CP}^2,\ \mathbb {CH}^2,\ S^2 \times S^2,\ S^2 \times \mathbb R^2,\ S^2 \times \mathbb H^2,\ \mathbb R^2 \times \mathbb H^2,\\
\mathbb H^2 \times \mathbb H^2,\
\til{SL}_2 \times \mathbb R,\ Nil^3 \times \mathbb R,\ Sol_0^4,\ Sol_1^4.
\end{align*}
\item $\mathbb R^4$ admits a complex structure compatible with $\mathbb R^4 \ltimes U_2$, and $S^3 \times \mathbb R$ admits a complex structure compatible with $U_2 \times \mathbb R$.  All other geometries admit no compatible complex structure.
\item In every case except $Sol_1^4$ the complex structure is unique up to isomorphism, and on $Sol_1^4$ there are two isomorphism classes.  As complex manifolds these are denoted $Sol_1^4$ and $(Sol_1^4)'$
\item There are K\"ahler metrics compatible with the complex and geometric structure precisely in the cases
\begin{align*}
\mathbb {CP}^2,\ \mathbb {CH}^2,\ S^2 \times S^2,\ S^2 \times \mathbb R^2,\ S^2 \times \mathbb H^2,\ \mathbb R^2 \times \mathbb H^2,\mathbb H^2 \times \mathbb H^2.
\end{align*}
\end{enumerate}
\end{thm}

\subsection{Existence and convergence results}

For locally homogeneous pluriclosed metrics, the pluriclosed flow will reduce to a system of ODE which greatly simplifies the analysis.  We record theorem statements here assuming familiarity with the Kodaira classification of surfaces, described in more detail in part II.  First there is a complete description of the long time existence behavior as well as the limiting behavior on the universal cover.

\begin{thm} \label{t:Boling1} (\cite{Boling} Theorem 1.1) Let $g_t$ be a locally homogeneous solution of pluriclosed flow on a compact complex surface which exists on a maximal time interval $[0,T)$.  If $T < \infty$ then the complex surface is rational or ruled.  If $T = \infty$ and the manifold is a Hopf surface, the evolving metric converges exponentially fast to a multiple of the Hopf metric (cf. \ref{f:Hopfmetric}).  Otherwise, there is a blowdown limit
\begin{align*}
\left(\til{g}_{\infty}\right)_t = \lim_{s \to \infty} s^{-1} \til{g}_{st}
\end{align*}
of the induced flow on the universal cover which is an expanding soliton in the sense that $\til{g}_t = t \til{g}_1$.
\end{thm}

It is implicit in the proof that one can actually recover the underlying complex surface from the asymptotic behavior of the flow.  In particular, if a locally homogeneous solution to pluriclosed flow encounters a finite time singularity, then the singularity is of type I and the underlying complex surface is rational or ruled.  If it is type IIb, i.e. $\brs{\Rm} t \to \infty$, then the underlying manifold is a diagonal Hopf surface, whereas if it is type III, i.e. $\brs{\Rm} t < \infty$, then the underlying manifold is one of the remaining surfaces in the Wall classification, i.e. a torus, hyperlliptic, Kodaira, or Inoue surface.  In fact, much sharper statements can be made concerning the Gromov-Hausdorff limits at infinity which recover precisely the underlying complex surface.  

\begin{thm} \label{t:Boling2} (\cite{Boling} Theorem 1.2) Let $\gw_t$ be a locally homogeneous solution of pluriclosed flow on a compact complex surface $(M, J)$ which exists on $[0,\infty)$ and suppose $(M, J)$ is not a Hopf surface.  Let $\hat{\gw}_t = \frac{\gw_t}{t}$.
\begin{enumerate}
\item If the surface is a torus, hyperelliptic, or Kodaira surface, then the family $(M, \hat{\gw}_t)$ converges as $t \to \infty$ to a point in the Gromov-Hausdorff topology.
\item If the surface is an Inoue surface, then the family $(M, \hat{\gw}_t)$ converges as $t \to \infty$ to a circle in the Gromov-Hausdorff topology and moreover the length of this circle depends only on the complex structure of the surface.
\item If the surface is a properly elliptic surface where the genus of the base curve is at least $2$, then the family $(M, \hat{\gw}_t)$ converges as $t \to \infty$ to the base curve with constant curvature metric in the Gromov-Hausdorff topology.
\item If the surface is of general type, then the family $(M, \hat{\gw}_t)$ converges as $t \to \infty$ to a product of K\"ahler--Einstein metrics on $M$.
\end{enumerate}
\end{thm}

\vskip 0.2in
\begin{center}
\textbf{\large{Part II: Geometrization of complex surfaces}}
\end{center}
\vskip 0.2in

First we recall the rudimentary aspects of the Kodaira classification of surfaces, referring the reader to the classic text \cite{BHPV} for a very detailed accounting.  The first main tool for classifying complex surfaces is the Kodaira dimension.  If $(M^4, J)$ is a complex surface, let $K$ denote the canonical bundle, and let $p_n = \dim H^0(K^{\otimes n})$ denote the plurigenera of $M$.  If all $p_n$ are zero, we say that the Kodaira dimension of $M$ is $\kod(M) = -\infty$.  Otherwise $\kod(M)$ is the smallest integer $k$ such that $\frac{p_n}{n^k}$ is bounded, which is no greater than $2$ in this case, giving the four possibilities, $-\infty, 0, 1$ and $2$.  Within these four classes one wants to understand which surfaces admit K\"ahler metrics.  Applying the K\"ahler identities to Dolbeault cohomology, one obtains that the odd Betti numbers of a K\"ahler surface are even.  Building on many results in the Kodaira classification, Siu \cite{SiuK3} proved the converse: if $b_1(M)$ is even then the surface admits K\"ahler metrics.  Later, direct proofs of this equivalence were given by Buchdahl \cite{BuchdahlK}, Lamari \cite{Lamari}.  This leaves in principle eight classes of surfaces, although it follows from Grauert's ampleness criterion (cf. \cite{BHPV} IV.6) that any surface with $\kod(M) = 2$ is projective, and hence K\"ahler, thus non-K\"ahler surfaces only occur for $\kod(M) \leq 1$, leaving seven classes.  We will address each of these classes in turn, and state the conjectural behavior of pluriclosed flow for each case.

\section{Conjectural limiting behavior on K\"ahler surfaces} \label{s:conjlimitsK}

It is natural to ask why one would bother studying pluriclosed flow on a K\"ahler surface, since one can detect via topological invariants whether a given complex surface admits K\"ahler metrics, and thus simply study K\"ahler-Ricci flow on these manifolds to produce canonical geometries.  However, as we ultimately want to apply pluriclosed flow as an a priori geometrizing process relying on minimal hypotheses, understanding its properties with arbitrary initial data even on K\"ahler manifolds plays an important philosophical role.  Moreover, K\"ahler surfaces are a rich class of complex manifolds on which to test the naturality and tameness of pluriclosed flow from a PDE point of view.  Most importantly, as we will see in part III, there are concrete applications to understanding the classification of generalized K\"ahler manifolds that cannot be approached through the K\"ahler-Ricci flow.

Having said this, the guiding principle here is quite simple: 
\begin{quote}
``\emph{Pluriclosed flow behaves like K\"ahler-Ricci flow on K\"ahler manifolds.}''
\end{quote}
To illustrate, note that Theorems \ref{t:negcurvthm} and \ref{t:torustheorem} yield the same existence time, and converge to the same limits as, the K\"ahler-Ricci flow in those settings.  The conjectures to follow below all follow this principle.

\begin{rmk} In the discussion below we make the assumption that the underlying complex surface is \emph{minimal}, i.e. free of $(-1)$-curves.  By a standard argument (\cite{BHPV} Theorem III.4.5) one can perform a finite sequence of blowdowns to obtain a minimal complex surface.  Thus from the point of view of complex geometry little is lost by considering only minimal surfaces.  Of course from the point of view of analysis one still would like to know what happens to the flow in the general setting.  An elementary calculation using the adjunction formula shows that pluriclosed flow homothetically shrinks the area of $(-1)$ curves to zero in finite time.  Conjecturally one expects that pluriclosed flow ``performs the blowdown'' in an appropriate sense.  This has been confirmed in some special cases for K\"ahler-Ricci flow \cite{SW}.
\end{rmk}

\subsection{Surfaces of general type}

By definition these are complex surfaces with Kodaira dimension $2$.  These surfaces have $c_1^2(M) > 0$ and so it follows from Grauert's ampleness criterion (cf. \cite{BHPV} IV.6) that all such are automatically projective, hence K\"ahler.  In the case $c_1 < 0$ the existence of a K\"ahler-Einstein metric follows from the work of Aubin-Yau, and the construction of this metric via K\"ahler-Ricci flow follows from Cao \cite{CaoKRF}.

\begin{thm} \label{t:Aubin-Yau-Cao} (Aubin-Yau \cite{Aubin} \cite{YauCC}, Cao \cite{CaoKRF}) Let $(M^4, J)$ be a compact K\"ahler surface with $c_1 < 0$.  There exists a unique metric $\gw_{\tiny{\mbox{KE}}} \in - c_1$ satisfying $\rho_{\gw_{\tiny{\mbox{KE}}}} = - \gw_{\tiny{\mbox{KE}}}$.  Moreover, given any K\"ahler metric $\gw_0$, the solution to K\"ahler-Ricci flow with initial condition $\gw_0$ exists on $[0,\infty)$, and
\begin{align*}
\lim_{t \to \infty} \frac{\gw_t}{t} = \gw_{\tiny{\mbox{KE}}}.
\end{align*}
\end{thm}
\noindent More generally, these surfaces can have $c_1 \leq 0$, admitting finitely many $-2$ curves, and their canonical models are orbifolds obtained by contraction of these curves.  Tian-Zhang \cite{TianZhang} proved in this setting that the K\"ahler-Ricci flow exists globally, converging after normalization to the unique oribfold K\"ahler-Einstein metric on the canonical model.  

To describe the conjectured behavior of the pluriclosed flow, we first compute the formal existence time.  By choosing a background metric with $\rho(\til{\gw}) \leq 0$, given any pluriclosed metric $\gw_0$ it follows that $\gw_0 - t \rho(\til{\gw}) > 0$, and so $\tau^*(\gw_0) = \infty$.  Thus we expect the pluriclosed flow to exist globally on such manifolds, and converge after normalization to the unique K\"ahler-Einstein metric on the canonical model.  Theorem \ref{t:negcurvthm} confirms this conjecture for a large class of surfaces of general type. 

\begin{conj} \label{c:gentypeconj} Let $(M, J)$ be a complex surface of general type.  Given $\gw_0$ a pluriclosed metric, the solution to pluriclosed flow with initial condition $\gw_0$ exists on $[0,\infty)$, and the solution to the normalized pluriclosed flow exists on $[0,\infty)$ and converges exponentially to $\gw_{KE}$, the unique orbifold K\"ahler-Einstein metric on the associated canonical model.
\end{conj}

\subsection{Properly Elliptic surfaces}

By definition these are K\"ahler complex surfaces of Kodaira dimension $1$.  For such surfaces there is a curve $\Sigma$ together with a holomorphic map $\pi : X \to \Sigma$ such that the canonical bundle of $M$ satisfies $K = \pi^* L$ for an ample line bundle $L$ over $\Sigma$.  Further, the generic fiber of $\pi$ is a smooth elliptic curve, and we call such points regular.  Near regular points we obtain a map to the moduli space of elliptic curves, and thus we can pull back the $L^2$ Weil-Petersson metric  to obtain a semipositive $(1,1)$ form $\gw_{WP}$ on the regular set.  This form plays a key role in describing the limiting behavior of K\"ahler-Ricci flow on these surfaces, done by Song-Tian\cite{SongTianPKD}:

\begin{thm} \label{t:SongTian} (\cite{SongTianPKD} Theorem 1.1) Let $\pi : M \to \Sigma$ be a minimal elliptic surface of $\kod(M) = 1$ with singular fibers $M_{s_1} = m_1 F_1, \dots, M_{s_k} = m_k F_k$ of multiplicity $m_i \in \mathbb N$, $i = 1\dots,k$.  Then for any initial K\"ahler metric $\gw_0$, the normalized K\"ahler-Ricci flow with this initial data exists on $[0,\infty)$ and satisfies:
\begin{enumerate}
\item $\gw_t$ converges to $\pi^* \gw_{\infty} \in -2 \pi c_1(M)$ as currents for a positive current $\gw_{\infty}$ on $\Sigma$.
\item $\gw_{\infty}$ is smooth on $\Sigma_{\mbox{\tiny{reg}}}$ and $\rho(\gw_{\infty})$ is a well-defined current on $\Sigma$ satisfying
\begin{align*}
\rho(\gw_{\infty}) = - \gw_{\infty} + \gw_{\mbox{\tiny{WP}}} + 2 \pi \sum_{i=1}^k \frac{m_k-1}{m_k} [s_i].
\end{align*}
\end{enumerate}
\end{thm}

Turning to pluriclosed flow, we first note that since $K_M$ is semiample, there exists a background metric with nonpositive Ricci curvature, and using this background metric it is clear that for any pluriclosed metric $\gw_0$ one has $\tau^*(\gw_0) = \infty$.  Furthermore, we conjecture the same limiting behavior as the K\"ahler-Ricci flow:

\begin{conj} \label{c:ellipticconj} Let $\pi : M \to \Sigma$ be a minimal elliptic surface of $\kod(M) = 1$.  Given $\gw_0$ a pluriclosed metric, the solution to normalized pluriclosed flow with initial condition $\gw_0$ exists on $[0,\infty)$, and satisfies conclusions $(1)$ and $(2)$ of Theorem \ref{t:SongTian}.
\end{conj}

\subsection{Elliptic surfaces of Kodaira dimension $0$}

It follows from the Kodaira classification that a minimal K\"ahler surface of Kodaira dimension $0$ is finitely covered by either a torus or a $K3$ surface.  In particular, $c_1 = 0$, and a Ricci flat (Calabi-Yau) metric is known to exist in each K\"ahler class by Yau \cite{YauCC}.  The global existence and convergence of K\"ahler-Ricci flow to a Calabi-Yau metric follows from Cao \cite{CaoKRF}:

\begin{thm} \label{t:Yau-Cao} (Yau \cite{YauCC}, Cao \cite{CaoKRF}) Let $(M^4, J)$ be a compact K\"ahler surface with $c_1 = 0$.  Given $[\gw]$ a K\"ahler class, there exists a unique metric $\gw_{\tiny{\mbox{CY}}} \in [\gw]$ satisfying $\rho_{\gw_{\tiny{\mbox{CY}}}} = 0$.  Moreover, given any K\"ahler metric $\gw_0$, the solution to K\"ahler-Ricci flow with initial condition $\gw_0$ exists on $[0,\infty)$, satisfies $[\gw_t] \in [\gw]$, and moreover
\begin{align*}
\lim_{t \to \infty} \gw_t = \gw_{\tiny{\mbox{CY}}}.
\end{align*}
\end{thm}

For the pluriclosed flow, first note that as $c_1 = 0$, the Aeppli cohomology class is fixed along pluriclosed flow, and so $\tau^*(\gw_0) = \infty$ for any $g_0$, thus following the general principle we conjecture global existence, as well as convergence to a Calabi-Yau metric.

\begin{conj} \label{c:ellipticconj2} Let $(M^4, J)$ be a compact K\"ahler surface with $\kod(M) = 0$.  Given $\gw_0$ a pluriclosed metric on $M$, the pluriclosed flow with initial condition $\gw_0$ exists on $[0,\infty)$, and converges to a K\"ahler Calabi-Yau metric $\gw_{CY}$.
\end{conj}

Theorem \ref{t:torustheorem} confirms this for tori, while (\cite{StreetsNDGKS} Theorem 1.1) establishes the long time existence and weak convergence for certain solutions to pluriclosed flow arising in generalized K\"ahler geometry (cf. \S \ref{s:GKRF} below).

\subsection{Rational and ruled surfaces}

By the Kodaira classification K\"ahler surfaces with $\kod(M) = -\infty$ are birational to $\mathbb{CP}^2$.  An important case of such surfaces are Fano surfaces, i.e. surfaces with $c_1 > 0$, as these are candidates to admit K\"ahler-Einstein metrics.  It follows from the Enriques-Kodaira classification that the smooth manifolds underlying Fano surfaces are $\mathbb {CP}^1 \times \mathbb{CP}^1$ and $\mathbb{CP}^2 \# k \bar{\mathbb{CP}}^2$, for $0 \leq k \leq 8$.  The existence of a K\"ahler-Einstein metric in this setting is obstructed in general, in particular Matsushima showed the automorphism group must be reductive (\cite{Matsushima}), with further obstructions due to Futaki \cite{Futaki}.  This rules out $\mathbb{CP}^2 \# k \bar{\mathbb{CP}}^2$ for $k=1,2$.  For the remaining cases, Tian settled the existence question:

\begin{thm} (\cite{TianCalabi}) Let $(M^4, J)$ be a compact Fano surface such that the Lie algebra of the automorphism group is reductive.  There exists a metric $\gw_{KE} \in c_1$ satisfying $\rho_{\gw_{KE}} = \gw_{KE}$.
\end{thm}

Extending our point of view slightly, we can look for K\"ahler-Ricci solitons, i.e. solutions of
\begin{align*}
\rho_{\gw} = \gw + L_{V} \gw,
\end{align*}
where $V$ is a holomorphic vector field.  Koiso \cite{Koisorot} showed the existence of a soliton on $\mathbb {CP}^2 \# \bar{\mathbb{CP}}^2$ using a symmetry ansatz.  Later Wang-Zhu \cite{WangZhuKRS} showed the existence of K\"ahler-Ricci solitons on all toric K\"ahler manifolds of positive first Chern class, in particular covering the case of $\mathbb {CP}^2 \# 2 \bar{\mathbb{CP}}^2$.

It is natural to expect the K\"ahler-Ricci flow to converge to these canonical metrics, when they exist.  There is a fairly complete picture of the K\"ahler-Ricci flow on Fano surfaces confirming this expectation:

\begin{thm} \label{t:FanoKRF} (\cite{TianZhuCOF,TZZZ}) Let $(M^4, J)$ be a compact Fano surface admitting a K\"ahler-Ricci soliton $(\gw_{KRS}, X)$.   Given $\gw_0 \in c_1$ a K\"ahler metric which is $G_X$-invariant, where $G_X$ is the one-parameter subgroup generated by $\mbox{Im}(X)$, the solution to normalized K\"ahler-Ricci flow exists on $[0,\infty)$ and converges to a K\"ahler-Ricci soliton.
\end{thm}

For the pluriclosed flow in the Fano setting, the normalized flow will fix the Aeppli cohomology class, and thus we expect global existence of the flow, and convergence to a soliton when one exists.

\begin{conj} \label{c:PCFFanoconj} Let $(M^4, J)$ be a compact Fano surface admitting a K\"ahler-Ricci soliton $(\gw_{KRS}, X)$.   Given $\gw_0 \in c_1$ a K\"ahler metric which is $G_X$-invariant, where $G_X$ is the one-parameter subgroup generated by $\mbox{Im}(X)$, the solution to normalized pluriclosed flow exists on $[0,\infty)$ and converges to a K\"ahler-Ricci soliton.
\end{conj}
\noindent As mentioned in the introduction, this conjecture would have consequences for the classification of generalized K\"ahler structures, and this will be detailed in part III below.  As we describe in Theorem \ref{t:CP2} below, the global existence of the normalized pluriclosed flow as well as a weak form of convergence can be established on $\mathbb{CP}^2$, which proves a kind of uniqueness for generalized K\"ahler structure on $\mathbb{CP}^2$.

Beyond Fano surfaces, the behavior of K\"ahler-Ricci flow more generally on rational and ruled surfaces can be quite delicate.  Even on a fixed Hirzebruch surface the K\"ahler-Ricci flow  can exhibit diverse behavior depending on the initial K\"ahler class.  In particular, Song-Weinkove (\cite{SW} Theorem 1.1) confirmed a conjecture of Feldman-Ilmanen-Knopf (\cite{FIK}) and showed that the flow can collapse $\mathbb{CP}^1$ fibers, shrink to a point, or contract exceptional divisors, depending on the choice of initial K\"ahler class.  More generally Song-Sz\'ekelyhidi-Weinkove (\cite{SSW}) established contraction of the fibers for more general ruled surfaces.  Rather than delving into these possibilities we refer to the general principle that one expects the same kinds of singularity formation for the pluriclosed flow.

\section{Conjectural limiting behavior on non-K\"ahler surfaces} \label{s:conjlimitsNK}

We now turn to the case of non-K\"ahler surfaces.  As it turns out, most complex surfaces outside of Class $\mbox{VII}_+$ admit locally homogeneous metrics, and thus admit special solutions to pluriclosed flow as described in Theorems \ref{t:Boling1}, \ref{t:Boling2}.  It is natural to expect that the convergence behavior on such manifolds for arbitrary initial data is the same as that in the homogeneous case, and we formalize and expand on this below.  Turning to the case of Class $\mbox{VII}_+$, we do not have homogeneous solutions available to guide us.  Nonetheless based on some examples and carefully inspecting monotonicity formulas for pluriclosed flow, we motivate a conjectural description of all long time solutions, which is loosely related to the theory of complete Calabi-Yau metrics on the complement of canonical divisors in surfaces, pioneered by Tian-Yau \cite{TianYau1, TianYau2}.

\subsection{Properly elliptic surfaces}

For elliptic fibrations with $b_1$ odd and $\kod(M) \geq 0$, it is known (cf. \cite{Wall} Lemma 7.2, \cite{Brinz2} Lemmas 1, 2) that only multiple fibers can occur, and thus there is a finite covering of $M$ which is an elliptic fiber bundle over a new curve which is a branched cover of the original base curve.  Furthermore, the condition $\kod(M) = 1$ is equivalent to the Euler characteristic of the base curve being negative, and these are called properly elliptic surfaces.  For these surfaces the canonical bundle is numerically effective, and so it follows from Corollary \ref{c:conecor} that $\tau^*(\gw_0) = \infty$ for any initial metric, and thus we expect global existence of pluriclosed flow.  Due to the lack of singular fibers, the expected limiting behavior in this case is actually simpler than that observed previously for the K\"ahler-Ricci flow on elliptic surfaces of Kodaira dimension $1$ (cf. Theorem \ref{t:SongTian}).  Also, for these surfaces it follows from (\cite{Wall} Theorem 7.4) that $M$ admits a geometric structure modeled on $\til{SL}_2 \times \mathbb R$.  Thus in particular Theorem \ref{t:Boling2} (3) gives a class of metrics exhibiting global existence and Gromov-Hausdorff convergence to the base curve.  Moreover, by Theorem \ref{t:Boling1}, there exists a blowdown limit of the solutions at infinity.  These exhibit the interesting behavior that the blowdown limits on the universal covers are invariant metrics on $\mathbb H \times \mathbb C$, so that the geometric type ``jumps'' in the limit.  This behavior is conjecturally universal.

\begin{conj} \label{c:propellipconj} Let $(M^4, J)$ be a properly elliptic surface with $\kod(M) = 1$ and $b_1$ odd.  Given $\gw$ a pluriclosed metric on $M$, the solution to pluriclosed flow with this initial data exists on $[0,\infty)$, and $(M, \frac{\gw_t}{t})$ converges in the Gromov-Hausdorff topology to $(B, \gw_{\mbox{\tiny{KE}}})$, the base curve with canonical orbifold K\"ahler-Einstein metric.  On the universal cover of $M$ there is a blowdown limit $\til{\gw}_{\infty}(t) = \lim_{s \to \infty} s^{-1} \til{\gw}(st)$ which is a locally homogeneous expanding soliton.
\end{conj}

Recently we were able to establish some cases of this conjecture, namely for initial data invariant under the torus action.

\begin{thm} (\cite{SPCFERYM}) Let $(M^4, J)$ be a properly elliptic surface with $\kod(M) = 1$ and $b_1$ odd.  Given $\gw$ a $T^2$-invariant pluriclosed metric on $M$, the solution to pluriclosed flow with this initial data exists on $[0,\infty)$, and $(M, \frac{\gw_t}{t})$ converges in the Gromov-Hausdorff topology to $(B, \gw_{\mbox{\tiny{KE}}})$, the base curve with canonical orbifold K\"ahler-Einstein metric.  On the universal cover of $M$ there is a blowdown limit $\til{\gw}_{\infty}(t) = \lim_{s \to \infty} s^{-1} \til{\gw}(st)$ which is a locally homogeneous expanding soliton.
\end{thm}

\subsection{Kodaira surfaces}

Kodaira surfaces are the complex surfaces with $\kod(M) = 0$ and $b_1$ odd.  These are elliptic fiber bundles over elliptic curves, called primary Kodaira surfaces, or finite quotients of such, called secondary Kodaira surfaces (cf. \cite{BHPV} V.5).  Wall showed (\cite{Wall} Theorem 7.4) that these surfaces admit a geometric structure modelled on $\Nil^3 \times \mathbb R$.  This was refined by Klingler \cite{Klingler}, who gave the precise description of the universal covers as well as the presentations of the lattice subgroups.  To describe the pluriclosed flow on these surfaces we first compute the formal existence time.  Topological considerations (cf. \cite{BHPV} pg. 197) show that the canonical bundle of primary Kodaira surfaces is trivial.  Thus $K \cdot D = 0$ for all divisors $D$, and thus for any pluriclosed metric $\gw_0$ it follows directly from Corollary \ref{c:conecor} that $\tau^*(\gw_0) = \infty$.  Thus we expect global existence of the flow in general.  Theorem \ref{t:Boling2} confirms this, and gives the limiting behavior of the flow, in the locally homogeneous setting, which we conjecture to be the general behavior.

\begin{conj} \label{c:Kodairaconj} Let $(M^4, J)$ be a Kodaira surface.  Given $\gw_0$ a pluriclosed metric on $M$, the solution to pluriclosed flow with this initial data exists on $[0,\infty)$, and $(M, \frac{\gw_t}{t})$ converges in the Gromov-Hausdorff topology to a point.  On the universal cover of $M$ there is a blowdown limit $\til{\gw}_{\infty}(t) = \lim_{s \to \infty} s^{-1} \til{\gw}(st)$ which is a locally homogeneous expanding soliton.
\end{conj}

Similar to the case of properly elliptic surfaces, we can exploit the fibration structure to prove the conjecture for invariant initial data.

\begin{thm} (\cite{SPCFERYM}) Let $(M^4, J)$ be a Kodaira surface.  Given $\gw$ a $T^2$-invariant pluriclosed metric on $M$, the solution to pluriclosed flow with this initial data exists on $[0,\infty)$, and $(M, \frac{\gw_t}{t})$ converges in the Gromov-Hausdorff topology to a point.  On the universal cover of $M$ there is a blowdown limit $\til{\gw}_{\infty}(t) = \lim_{s \to \infty} s^{-1} \til{\gw}(st)$ which is a locally homogeneous expanding soliton.
\end{thm}

\subsection{Class $\mbox{VII}_0$ surfaces}

Surfaces of Class $\mbox{VII}_0$ are defined by the conditions $\kod(M) = -\infty$, $b_1 = 1$, and $b_2 = 0$.  Hopf surfaces are of this type, and Inoue \cite{Inouesurfaces} gave classes of examples, and showed that surfaces in this class with no curves must be one of his examples.  Later, using the methods of gauge theory, Li-Yau-Zheng \cite{LiYauBog} and Teleman \cite{TelemanBog} showed that all surfaces in this class are either Hopf surfaces or Inoue surfaces.  We describe the conjectured behavior on these two classes of surfaces below.

\subsubsection{Inoue surfaces}

Inoue surfaces were introduced in \cite{Inouesurfaces}, giving examples of complex surfaces with $b_1 = 1, b_2 = 0$, and containing no curves.  These come in three classes, all of which are quotients of $\mathbb H \times \mathbb C$ by affine subgroups, and we will recall the construction of the simplest of these classes.  First, fix $Z \in \SL(3,\mathbb Z)$, with eigenvalues $\ga, \gb, \bgb$ such that $\ga > 1$ and $\gb \neq \bgb$.  Choose a real eigenvector $(a_1,a_2,a_3)$ for $\ga$ and an eigenvector $(b_1,b_2,b_3)$ for $\gb$.  It follows that the three complex vectors $(a_i, b_i)$ are linearly independent over $\mathbb R$.  Using these we define a group of automorphisms $G_Z$ of $\mathbb H \times \mathbb C$ generated by
\begin{align*}
g_0 (w,z) =&\ (\ga w, \gb z),\\
g_i(w,z) =&\ (w + a_i, z + b_i), \qquad i = 1,2,3.
\end{align*}
This action is free and properly discontinuous, and so defines a quotient surface $S_Z = \mathbb H \times \mathbb C/G_Z$. If we let $G$ denote the subgroup generated by $g_1,g_2,g_3$, it is clear that this is isomorphic to $\mathbb Z^3$, and preserves the affine varieties
\begin{align*}
(w_0,z_0) + (\mathbb R (a_1,b_1) \oplus \mathbb R(a_2,b_2) \oplus \mathbb R(a_3,b_3)).
\end{align*}
These spaces are parameterized by $\Im(w_0)$, and thus the quotient $\mathbb H \times \mathbb C / G$ is diffeomorphic to $T^2 \times \mathbb R_+$.  Since $g_0$ preserves the fibers of this quotient, it follows that the quotient $S_Z$ is a three-torus bundle over $S^1$, with $b_1 = 1$, $b_2 = 0$.  Wall shows (\cite{Wall} Proposition 9.1) that these surfaces are precisely those with geometric structure $Sol_0^4$, $Sol_1^4$, or $(Sol_1^4)'$.  In particular for the example above one sees that each $g_i \in Sol_0^4$.

To understand the pluriclosed flow on these surfaces, we first compute the formal existence time.  As shown by Inoue (\cite{Inouesurfaces} \S 2-\S 4), these surfaces contain no curves, thus it follows immediately from Corollary \ref{c:conecor} that $\tau^*(\gw_0) = \infty$ for any $\gw_0$.  This long time existence was verified in the homogeneous setting by Boling (Theorem \ref{t:Boling1}), and also in the commuting generalized K\"ahler setting (\cite{StreetsCGKFlow} Theorem 1.3, cf. Theorem \ref{t:CGKRF} below).  Boling also established interesting convergence behavior in the homogeneous setting (Theorem \ref{t:Boling2}), showing that the pluriclosed flow collapses the three-torus fibers described above, yielding the base circle as the Gromov-Hausdorff limit.  Moreover, he showed the existence of a blowdown limit on the universal cover which is an expanding soliton.  We conjecture that this is the general behavior.

\begin{conj} \label{c:Inoueconj} Let $(M^4, J)$ be an Inoue surface.  Given $\gw_0$ a pluriclosed metric on $M$, the solution $\gw_t$ to pluriclosed flow with this initial data exists on $[0,\infty)$.  The family $(M, \frac{\hat{\gw}_t}{t})$ converges as $t \to \infty$ to a circle in the Gromov-Hausdorff topology and moreover the length of this circle depends only on the complex structure of the surface.  On the universal cover of $M$ there is a blowdown limit $\til{\gw}_{\infty}(t) = \lim_{s \to \infty} s^{-1} \til{\gw}(st)$ which is a canonical locally homogeneous expanding soliton.
\end{conj}

\subsubsection{Hopf surfaces} \label{ss:Hopf}

Hopf surfaces by definition are all compact quotients of $\mathbb C^2 \backslash \{0\}$.  A Hopf surface is called primary if its fundamental group is isomorphic to $\mathbb Z$.  Since every Hopf surface has a finite cover which is primary, we restrict our discussion to the primary case.  By (\cite{Kod2} Theorem 30), the possible group actions on $\mathbb C^2 \backslash \{0\}$ take the form
\begin{align*}
\gg(z_1,z_2) =&\ (\ga z_1, \gb z_2 + \gl z_1^m),\ \mbox{ where } 0 < \brs{\ga} \leq \brs{\gb} < 1, \qquad (\ga - \gb^m)\gl = 0.
\end{align*}
We call the resulting quotient surface $(M, J_{\ga,\gb,\gl})$.  It follows that $M$ is always diffeomorphic to $S^3 \times S^1$, thus $b_1=1, b_2 = 0$, and the surfaces are not K\"ahler.  We call the surface Class $1$ if $\gl = 0$, and Class $0$ otherwise.  Within Class $1$ Hopf surfaces we call the surface diagonal if $\brs{\ga} = \brs{\gb}$.  

To describe the conjectured behavior of pluriclosed flow, we begin with a simple example.  On $\mathbb C^2 \backslash \{0\}$, define the metric
\begin{align} \label{f:Hopfmetric}
\gw_{\mbox{\tiny{Hopf}}} := \frac{\i \del \delb \brs{z}^2}{\brs{z}^2}.
\end{align}
Written in this form, it is clear that this is a Hermitian metric which is invariant under $\gg(z_1,z_2) = (\ga z_1, \gb z_2)$, where $\brs{\ga} = \brs{\gb}$, and thus descends to a metric on the diagonal Hopf surfaces.  Elementary calculations show that this metric is pluriclosed.  Of course $\mathbb C^2 \backslash \{0\} \cong S^3 \times \mathbb R$, and as it turns out the Riemannian metric corresponding to $\gw_{\mbox{\tiny{Hopf}}}$ is isometric to $g_{S^3} \oplus ds^2$, where $g_{S^3}$ is the round metric on $S^3$ and $s$ is a parameter on $\mathbb R$.  Further calculations reveal that the torsion tensor $H$ is a multiple of $dV_{g_{S^3}}$, the volume form on $S^3$.  Since $\theta = \star H$ for complex surfaces, it follows that $\theta$ is a multiple of $ds$, and is in particular parallel.  Putting these facts together and applying Theorem \ref{t:PCFRGF}, it follows that $\gw_{\mbox{\tiny{Hopf}}}$ is a fixed point of pluriclosed flow, and in particular the pair $(g,H)$ satisfies
\begin{align*}
 \Rc - \tfrac{1}{4} H^2 =&\ 0,\\
d^* H =&\ 0.
\end{align*}
Thus we see explicitly a basic principle which appears to distinguish the non-K\"ahler setting from the K\"ahler setting: for this metric the torsion acts as a ``balancing force'' which cancels out the positive Ricci curvature of $g_{S^3}$.  Thus this metric is Bismut-Ricci flat, while certainly not Ricci flat.  In fact more is true: the Bismut connection is flat.  As $S^3$ is a simple Lie group, there are flat connections $\N^{\pm}$ which are compatible with a bi-invariant metric and having torsion $T(X,Y) = \pm [X,Y]$.  The connection $\N^+$, after taking a direct sum with a flat connection on $S^1$, recovers the Bismut connection of $\gw_{\mbox{\tiny{Hopf}}}$.  It follows from work of Cartan-Schouten (\cite{CartanSchouten, 
CartanSchouten2}), that triples $(M^{n}, g, H)$ with flat Bismut connection are isometric to products of simple Lie groups and classically flat spaces, as in this example.

Surprisingly, the metrics described above, which exist on diagonal Hopf surfaces, are the \emph{only} non-K\"ahler fixed points of pluriclosed flow, a result of Gauduchon-Ivanov (\cite{GauduchonIvanov} Theorem 2).  The central point is to employ a Bochner argument to show that the Lee form is parallel.  Thus it either vanishes identically, in which case the metric is K\"ahler and Calabi-Yau, or it is everywhere nonvanishing and yields a metric splitting of the universal cover.  Looking back at the fixed point equation, one then observes that the geometry transverse to the Lee form has constant positive Ricci curvature, and so is a quotient of the round sphere.  Thus the universal cover is isometric to $g_{S^3} \oplus ds^2$.  This identifies the metric structure, and using a lemma of Gauduchon (\cite{GauduchonWeyl} III Lemma 11) one identifies the possible complex structures, yielding only the diagonal Hopf surfaces.

This rigidity result is surprising since it implies that a naive extension of the Calabi-Yau theorem/Cao's theorem is impossible in the context of pluriclosed flow.  In particular, all Hopf surfaces have $b_2 = 0$ and thus $c_1 = 0$, and so $\tau^*(\gw_0) = \infty$ for all $\gw_0$, and so one expects global existence of the pluriclosed flow.  However, only the restricted class of diagonal Hopf surfaces admits fixed points.  This observation, together with the Perelman-type monotonicity formula, inspired the author to look for soliton-type fixed points of pluriclosed flow, as described in \S \ref{s:PCFsec}.  In \cite{Streetssolitons} we were able to construct such soliton solutions of pluriclosed flow on all Class $1$ Hopf surfaces, and show that Hopf surfaces are the only compact complex surfaces on which solitons can exist.  Note that this is a key qualitative distinction from the case of Ricci flow, where all compact steady solitons are trivial, i.e. Einstein.  While we were not able to rule out existence on Class $0$ Hopf surfaces, we conjecture that they do not exist.  It would be interesting to develop invariants akin to Futaki's invariants for K\"ahler-Einstein metrics to try to rule them out, or for that matter to rule out fixed points on non-diagonal Class $1$ Hopf surfaces via such invariants.

With this background we can now state the conjectured convergence behavior.  In particular, on diagonal Hopf surfaces the flow should converge to $\gw_{\mbox{\tiny{Hopf}}}$.  More generally, for Class $1$ Hopf surfaces one expects convergence to a soliton, which is presumably unique.  If it is true that solitons cannot exist on Class $0$ Hopf surfaces, one cannot obtain convergence to a soliton in the usual sense.  Due to the nature of Cheeger-Gromov convergence, it is however possible for the complex structure to ``jump'' in the limit.  In particular, every Class $1$ Hopf surface occurs as the central fiber in a family of Hopf surfaces, all other fibers of which are biholomorphic to the same Class $0$ Hopf surface.  The diffeomorphism actions necessary in taking Cheeger-Gromov limits on Class $0$ surfaces should result in the limiting complex structure jumping to the central fiber, i.e. the associated Class $1$ Hopf surface.

\begin{conj} \label{c:Hconj} Let $(M^4, J)$ be a primary Hopf surface.  Given $\gw_0$ a pluriclosed metric on $M$, the solution $\gw_t$ to pluriclosed flow exists on $[0,\infty)$.
\begin{enumerate}
\item If $(M, J)$ is a diagonal Hopf surface, $\gw_t$ converges in the $C^{\infty}$ topology to $\gw_{\mbox{\tiny{Hopf}}}$.
\item If $(M, J_{\ga,\gb})$ is a Class $1$ Hopf surface, $(M, J_{\ga,\gb}, \gw_t)$ converges in the $C^{\infty}$ Cheeger-Gromov topology to a unique steady soliton $(M, J_{\ga,\gb}, \gw_S)$.
\item If $(M, J_{\ga,\gb,\gl})$ is a Class $0$ Hopf surface, $(M, J_{\ga,\gb,\gl}, \gw_t)$ converges in the $C^{\infty}$ Cheeger-Gromov topology to a unique steady soliton $(M^4, {J}_{\ga,\gb}, \gw_S)$ on the Class $1$ Hopf surface adjacent to $J_{\ga,\gb,\gl}$ as described above.
\end{enumerate}
\end{conj}

Recently, in line with the results described above for properly elliptic and Kodaira surfaces, we are able to establish partial results confirming this behavior.  In this case, for technical reasons, to establish the global existence we require an extra condition, namely that the torsion is nowhere vanishing.

\begin{thm} (\cite{SPCFERYM}) Let $(M^4, J)$ be a diagonal Hopf surface.  Given $\gw_0$ a $T^2$-invariant pluriclosed metric with nowhere vanishing torsion, the solution to pluriclosed flow with this initial data exists on $[0,\infty)$ and converges to $\gw_{\mbox{\tiny{Hopf}}}$.
\end{thm}

\subsection{Class $\mbox{VII}_+$ surfaces}

Surfaces of Class $\mbox{VII}_+$ are defined by the conditions $\kod(M) = -\infty$, $b_1 = 1$, and $b_2 > 0$.  We briefly recall a general construction of surfaces of this type due to  Kato \cite{KatoGSS}, building on the initial construction of Inoue \cite{InoueNew}.  Let $\Pi_0$ denote blowup of the origin of the unit ball $B$ in $\mathbb C^2$.  Let $\Pi_1$ denote the blowup of a point $O_0 \in C_0 := \Pi_0^{-1}(0)$.  Iteratively let $\Pi_{i+1}$ denote blowup of a point $O_{i} \in C_i = \Pi_i^{-1}(O_{i-1})$.  Let $\Pi : \hat{B} \to B$ denote the composition of these blowups.  Choose a holomorphic embedding $\gs : \bar{B} \to \hat{B}$ such that $\gs(0) \in C_k$, the final exceptional divisor.  Let $N = \hat{B} \backslash \gs(\bar{B})$, which has two boundary components $\del \hat{B}$ and $\gs(\del B)$.  The map $\gs \circ \Pi : \del \hat{B} \to \gs(\del B)$ can be used to glue these two boundaries, producing a minimal compact complex surface $M = M_{\pi,\gs}$, and such surfaces are referred to as \emph{Kato surfaces}.  These surfaces are diffeomorphic to $S^3 \times S^1 \# k \mathbb{CP}^2$, but the complex structures are minimal.  It was conjectured by Nakamura (\cite{NakamuraClass} Conjecture 5.5) that \emph{all} complex surfaces of Class $\mbox{VII}_+$ are in fact Kato surfaces, and this remains the main open question in the Kodaira classification of surfaces. 

One approach to resolving this conjecture focuses on a particular geometric feature shared by all Kato surfaces, that of a \emph{global spherical shell} (GSS).  This is a biholomorphism of an annulus in $\mathbb {C}^2 - \{0\}$ into $M$ such that the image does not disconnect $M$.  These are easily seen to exist in Hopf surfaces, and moreover in the construction of Kato surfaces as above, any annulus around the origin in $B$ is a GSS.  The relevance of GSS was exhibited by Kato \cite{KatoGSS}, who showed that every surface admitting a global spherical shell is a degeneration of a blown up primary Hopf surface, and moreover is a Kato surface.  Later it was shown by Dloussky, Oeljeklaus, and Toma \cite{DOT} that if a complex surface of Class $\mbox{VII}_+$ admits $b_2$ complex curves, then it has a global spherical shell, and is hence a Kato surface.  Given this, Teleman \cite{TelemanDonaldson} proved the existence of a curve on all Class $\mbox{VII}_+$ surfaces with $b_2 = 1$, thus finishing their classification.  Moreover, when $b_2 = 2$ Teleman \cite{TelemanInstantons} showed the presence of a cycle of rational curves, again yielding the classification for this case.  This deep work remains the only proof of classification of these surfaces for $b_2 > 0$.

In our initial joint work with Tian \cite{PCF}, we conceived pluriclosed flow in part to address this classification problem.  In the follow-up \cite{PCFReg}, we discussed an argument by contradiction whereby a Class $\mbox{VII}_+$ with no curves at all would violate natural existence conjectures for the pluriclosed flow.  However, this line of argument via contradiction still leaves us with the question of what, even conjecturally, the flow may actually do in these settings.  First we address the existence time.  It follows from (\cite{Kod1} p. 755, \cite{Kod2} p. 683) that one has the following topological characteristics for any Class $\mbox{VII}_+$ surface:
\begin{align} \label{f:C7top}
h^{0,1} = 1, \quad h^{1,0} = h^{2,0} = h^{0,2} = 0, \quad b_2^+ = 0, \quad c_1^2 = - b_2.
\end{align}
Since the intersection form is negative definite, it follows from the adjunction formula that $K \cdot D \geq 0$ for any divisor $D$, with equality if and only if $D$ is a $(-2)$-curve.  Corollary \ref{c:conecor} then implies that $\tau^*(\gw_0) = \infty$ for any $\gw_0$ (cf. \cite{PCFReg} Proposition 5.7).

To describe the conjectural limiting behavior, we recall basic facts on the configuration of curves in these surfaces.  In the simplest case of parabolic Inoue surfaces, corresponding to a generic sequence of blowups in the construction described above, there is a cycle of $(-2)$ curves $C_i$ satisfying $K \cdot C_i = 0$, as well as an elliptic curve $E$ satisfying $K \cdot E = 1$.  In general, by blowing up the same point several times in the Kato construction, one can have rational curves of high negative self intersection, which then satisfy $K \cdot C > 0$.  We define
\begin{align*}
\Sigma_{0} := \left\{\mbox{curves } C,\ K \cdot C = 0 \right\},\\
\Sigma_{> 0} := \left\{\mbox{curves } C,\ K \cdot C > 0 \right\}.
\end{align*}
It follows from the adjunction formula that $\Sigma_0$ consists of smooth rational $(-2)$ curves, as well as rational curves with an ordinary double point and zero self-intersection.

Given these remarks on the structure of curves, a natural guess arises for the limiting behavior of pluriclosed flow which we now rule out.  Along pluriclosed flow the area of curves in $\Sigma_{0}$ will remain fixed, while the area of curves in $\Sigma_{>0}$ will grow linearly.  This is similar to the situation for surfaces of general type, where $(-2)$ curves remain fixed and all other curves grow linearly.  In that situation, one considers the normalized flow, for which the areas of $(-2)$ curves decay exponentially and the areas of all other curves approach a fixed positive value.  It was proved by Tian-Zhang (\cite{TianZhang}) that in this setting the K\"ahler-Ricci flow converges to a K\"ahler-Einstein metric on the canonical model of the original surface, which is the orbifold given by contraction of all $(-2)$ curves.  Thus one might expect a similar picture for pluriclosed flow, with the normalized flow converging to a canonical metric on an orbifold obtained by contraction of the $(-2)$ curves on the original surface.  However, it follows from a short calculation (cf. \cite{PCF} Proposition 3.8) that along pluriclosed flow the integral of the Chern scalar curvature evolves by $- c_1^2 = b_2 > 0$.  Thus for the normalized flow it approaches the value $b_2$.  However, using the expanding entropy functional for generalized Ricci flow (\cite{Streetsexpent}), if the normalized flow converged to some smooth metric on an orbifold, it would be a negative scalar curvature K\"ahler-Einstein metric, contradicting that the integral is $b_2 > 0$.

Instead, let us argue proceeding from the observation that the integral of Chern scalar curvature is asymptotically $b_2 t$.   As this implies that the scalar curvature is becoming postive on average, this would force the volume to go to zero if not for the torsion acting as a ``restoring force'' as described in \S \ref{ss:Hopf}.  Thus, for a point $p$ where the scalar curvature is positive, bounded and bounded away from zero, one expects to be able to construct a nonflat limit of pointed solutions $(M, \gw_t, J, p)$.  Arguing formally using the $\FF$-functional monotonicity (\ref{f:Fmon}) one expects this limit to be a steady soliton.  In almost every case, the set $\Sigma_{ > 0}$ is nonempty, and the areas of these curves grow linearly and so form the infinity of the resulting complete metric.  In these cases one expects that limits bases at any point $p \notin \Sigma_{> 0}$ will yield a nonflat steady soliton, while points in $\Sigma_{> 0}$ will yield a flat limit.  The exceptional cases are the Enoki surfaces described below which are certain exceptional compactifications of line bundles over elliptic curves.   In the special case of parabolic Inoue surfaces the zero section of the line bundle is an elliptic curve, and forms $\Sigma_{> 0}$.  Outside of this special case the set $\Sigma_{> 0}$ is empty, but nonetheless we expect some section of the line bundle to play the role of $\Sigma_{ > 0}$ as described above.  We summarize:

\begin{conj} \label{c:Class7conj} Let $(M^4, J)$ be a compact surface of Class $\mbox{VII}_+$.  Given $\gw_0$ a
pluriclosed metric on $M$ the solution $\gw_t$ to pluriclosed flow with this initial data exists on
$[0,\infty)$.  For a generic point $p \in M$ as described above, the pointed spaces $(M, J,\gw_t, p)$ converge in the pointed $C^{\infty}$ Cheeger-Gromov sense to a nonflat
complete steady soliton $(M_{\infty}, J_{\infty}, \gw_{\infty}, p)$ for some smooth function $f_{\infty}$.  Furthermore,
\begin{enumerate}
\item The vector field $\theta^{\sharp}_{\infty} + \N f_{\infty}$ is $J_{\infty}$-holomorphic.
\item $M_{\infty}$ admits a compactification to a complex surface $(\bar{M}_{\infty}, \bar{J}_{\infty})$.
\item The distribution orthogonal to $\theta^{\sharp}_{\infty} + \N f_{\infty}$ is integrable, and its generic leaf is an embedded submanifold of $M_{\infty}$, whose closure in $\bar{M}_{\infty}$ is a global spherical shell.
\end{enumerate}
\end{conj}

This is a natural extension of Conjecture \ref{c:Hconj}.  We note that for all steady solitons the vector field $\theta^{\sharp} + \N f$ is automatically holomorphic (\cite{Streetssolitons} Proposition 3.4, which can be extended to the complete setting).  Moreover, the twisted Lee form $e^{-f} (\theta + df)$ is automatically closed, and so this distribution is always integrable.  Also, for the solitons on Hopf surfaces constructed in \cite{Streetssolitons}, one can verify that the leaves of the distribution $\ker e^{-f} (\theta + df)$ are global spherical shells.  It would be interesting to determine if this was true for a general, non-K\"ahler complete steady soliton.  Given the nature of the convergence process, assuming the generic leaf is a global spherical shell, this will yield the existence of one on the original complex surface $(M, J)$, thus finishing the classification as described above.

Note that these conjectured limits bear a family resemblance to complete Calabi-Yau metrics, pioneered in work of Tian-Yau \cite{TianYau1, TianYau2}.  A prototypical result of this kind exhibits a complete Calabi-Yau metric on $M \backslash D$, where $M$ is smooth quasi-projective, and $D$ is a smooth anticanonical divisor with $D^2 \geq 0$.  Very loosely speaking the anticanonical divisor is a topological obstruction to the existence of a Ricci flat metric, and by removing it one can construct complete examples.  Similarly, in this setting the curves in $\Sigma_{> 0}$ form an obstruction to the existence of a ``Calabi-Yau'' type metric, which the flow naturally pushes to infinity.  We note however that complete Calabi-Yau do not arise as limits of the K\"ahler-Ricci flow on compact complex surfaces.  The difference in the expected qualitative behavior can be traced back to the simple fact that $c_1^2 = - b_2 \leq 0$ on these surfaces, which has a profound effect on the existence time and qualitative behavior on Class $\mbox{VII}_+$ surfaces versus other surfaces of Kodaira dimension $-\infty$.

For complete pluriclosed solitons as described, in line with a conjecture of Yau for complete Calabi-Yau metrics \cite{YauICM}, it is natural to expect that these complete solitons admit compactifications $(\bar{M}_{\infty}, \bar{J}_{\infty})$.  These should be compact complex surfaces with $b_1 = 1$, thus Class $\mbox{VII}_+$ surfaces, although here again it is possible that the complex structure has ``jumped'' in the limit so that this compactification need not be biholomorphic to the original surface.  Lastly, we note that not every Class VII surface admits holomorphic vector fields, and those that do are classified in \cite{DlousskyVF1, DlousskyVF2}.  Despite this restriction, the conjectured convergence can still happen since we do not necessarily expect the vector field to extend smoothly to the compactification, and moreover the complex structure can change in the limit as mentioned.

Let us flesh this picture out for specific classes of Kato surfaces.  First we consider the Enoki/Inoue surfaces, described as exceptional compactifications of line bundles over elliptic curves.  We let $E = \mathbb C^* / \IP{\ga}$, $0 < \brs{\ga} < 1$, be an elliptic curve.  Fix some $n \geq 1$ and $t \in \mathbb C^n$, identified with a polynomial via $t(w) = \sum t_k w^k$.  Using these one defines an automorphism of $\mathbb C \times \mathbb C^*$,
\begin{align*}
g_{n,\ga,t} (z,w) = (w^n z + t(w), \ga w).
\end{align*}
Let $A_{n,\ga,t}$ denote the quotient surface $\mathbb C \times \mathbb C^* / \IP{g_{n,\ga,t}}$, which is an affine line bundle over $E$.  Enoki showed \cite{Enoki}, generalizing a previous construction of Inoue \cite{InoueNew}, that these can be compactified with a cycle of rational curves, yielding a compact surface $S_{n,\ga,t}$ of Class $\mbox{VII}_+$.  Generically this cycle of curves are the only curves present in the surface, but in the case $t = 0$ (also known as parabolic Inoue surfaces), the zero section of the line bundle is a smooth elliptic curve in $S_{n,\ga,0}$, satisfying $K \cdot E = 1$.  Thus according to the conjecture above for parabolic Inoue surfaces we expect the area of this elliptic curve to grow linearly and ultimately form the ``infinity'' of a limiting complete steady soliton.  For more general Enoki surfaces $E_{n,\ga,t}$, there is still a topological torus in the homotopy class of the cycle of rational curves, whose area goes to infinity, conjecturally resulting in the same  steady soliton limit resulting from $E_{n,\ga,0}$.  It seems likely that the compactified limiting surfaces are all biholomorphic to $E_{n,\ga,0}$.

Next consider the Inoue-Hirzebruch surfaces.  These were initially constructed by Inoue \cite{InoueNew2}, using ideas related to Hirzebruch's description of Hilbert modular surfaces \cite{Hirz}.  These surfaces are constructed by resolving singularities of compactified quotients of $\mathbb H \times \mathbb C$.  These surfaces always admit two cycles of rational curves (cf. \cite{DlousskyIH} for a description of these surfaces and the structure of their curves), each of which contains curves of self intersection $\leq -3$.  Thus according to the conjecture above we expect the limit soliton to have at least two ends, possibly more depending on the structure of the curves.  In a given cycle of rational curves smooth $-2$ curves may intersect curves of high self intersection.  As their area stays fixed while the high self-intersection curve is pushed to infinity, this suggests that the $-2$ curve forms a finite area cusp end with the intersection point now at infinity.

The two examples of Enoki surfaces and Inoue-Hirzebruch surfaces form the extreme cases of Dloussky's index invariant \cite{DlousskyKato}, with the remaining cases called ``intermediate surfaces''.  Despite the intricate structure of curves on these surfaces, the limiting picture is essentially the same as described in the examples above, with the curves in $\Sigma_{> 0}$ being pushed to infinity, forming potentially several ends, and all $-2$ curves either contained in a compact region or intersecting curves at infinity with a cusp end.  

\vskip 0.2in
\begin{center}
\textbf{\large{Part III: Classification of generalized K\"ahler structures}}
\end{center}
\vskip 0.2in

\section{Generalized K\"ahler Geometry} \label{s:GKG}

Generalized K\"ahler geometry arose in work of Gates-Hull-Ro\v cek \cite{GHR}, in the course of their investigation of supersymmetric sigma models.  Later, using the framework of Hitchin's generalized geometry \cite{HitchinGCY}, Gualtieri \cite{GualtieriGCG, GualtieriGKG} understood generalized K\"ahler geometry in terms of a pair of commuting complex structures on $TM \oplus T^*M$ satisfying further compatibility conditions.  For our purposes here we will restrict ourselves to the ``classical'' formulation and not exploit the language of generalized geometry.  Thus a generalized K\"ahler structure on a manifold $M$ is a triple $(g, I, J)$ of a Riemannian metric together with two integrable complex structures, such that
\begin{align*}
d^c_I \gw_I = - d^c_J \gw_J, \qquad d d^c_I \gw_I = 0,
\end{align*}
where $\gw_I = g(I \cdot, \cdot)$, and $d^c_I = \i (\delb_I - \del_I)$, with analogous definitions for $J$.  

Associated to every generalized K\"ahler structure is a Poisson structure
\begin{align} \label{f:sigmadef}
\gs = \tfrac{1}{2} g^{-1} [I,J].
\end{align}
As shown by Pontecorvo \cite{PontecorvoCS} and Hitchin \cite{HitchinPoisson}, $\gs$ is the real part of a holomorphic Poisson structure with respect to both $I$ and $J$, in other words
\begin{align*}
\delb_I \gs^{2,0}_I = 0, \qquad \delb_J \gs^{2,0}_J = 0.
\end{align*}

\noindent The vanishing locus of $\gs$ has profound implications for the structure of generalized K\"ahler manifolds, and it is natural to understand their classification in terms of its structure.  The simplest case occurs when $\gs \equiv 0$.  In this case $[I,J] = 0$, and the endomorphism $Q = IJ$ satisfies $Q^2 = \Id$.  Thus $Q$ has eigenvalues $\pm 1$, and the eigenspaces of $Q$ yield a splitting $TM = T_+ \oplus T_-$.  These summands are $I$-invariant, thus we obtain a further splitting
\begin{align*}
T_{\mathbb C} M = T^{1,0}_+ \oplus T^{0,1}_+ \oplus T^{1,0}_- \oplus T^{0,1}_-.
\end{align*}
Complex surfaces admitting a holomorphic splitting of the tangent bundle were classified by Beauville \cite{BeauvilleSplit}.  Apostolov and Gualitieri \cite{ApostolovGualtieri} determined precisely which of these admits generalized K\"ahler structure of this type.

\begin{thm} \label{t:AG} (\cite{ApostolovGualtieri} Theorem 1) A compact complex surface $(M, I)$ admits a generalized K\"ahler structure $(g,I,J)$ with $[I,J] = 0$ if and only if $(M, I)$ is biholomorphic to
\begin{enumerate}
\item A ruled surface which is the projectivization of a projectively flat holomorphic vector bundle over a compact Riemann surface,
\item A bi-elliptic surface,
\item A surface of Kodaira dimension $1$ with $b_1$ even, which is an elliptic fibration over a compact Riemann surface, with singular fibers only multiple smooth elliptic curves,
\item A surface of general type, whose universal cover is biholomorphic to $\mathbb H \times \mathbb H$, with fundamental group acting diagonally on the factors.
\item A class $1$ Hopf surface
\item An Inoue surface of type $S_M$.
\end{enumerate}
\end{thm}

The next simplest case occurs when $\gs$ defines a nondegenerate bilinear form at all points.  It follows that $\Omega = \gs^{-1}$ is a symplectic form, which is the real part of a holomorphic symplectic form with respect to both $I$ and $J$.  The simplest example comes from hyperK\"ahler geometry.  If $(M^{4n}, g, I, J, K)$ is hyperK\"ahler, then the triple $(M^{4n}, g, I, J)$ is generalized K\"ahler, and using the quaternion relations one computes that $\Omega = \gw_K$.  As we will describe in \S \ref{ss:ndgk}, it is possible to deform this example to obtain a non-K\"ahler generalized K\"ahler structure with nondegenerate Poisson structure.  The existence of a holomorphic symplectic form places strong restritions on the underlying complex manifolds.  For complex surfaces, it follows from (\cite{AGG} Proposition 2) that the only possible underlying complex surfaces are tori, $K3$ surfaces, or primary Kodaira surfaces.  Considering the $I$ and $J$-imaginary pieces of $\Omega$, one overall obtains three independent self-dual forms, ruling out Kodaira surfaces which have $b_2^+ = 2$.  Thus such structures exist only on tori and $K3$ surfaces.

Generally, the Poisson structure can experience ``type change'', that is, the rank can drop on some locus.  As we are in four dimensions, and the rank jumps in multiples of $4$, the Poisson structure will be nondegenerate outside of a locus 
\begin{align*}
T = \left\{ p \in M\ |\ I = \pm J \right\}.
\end{align*}
It turns out that $T$ is a complex curve in both $(M, I)$ and $(M, J)$.  In the case the underlying surfaces are K\"ahler, $T$ is the support of an anticanonical divisor, and so only Del Pezzo surfaces are possible backgrounds, and Hitchin \cite{HitchindelPezzo} constructed generalized K\"ahler structures on these surfaces.  If the underlying surfaces are non-K\"ahler $T$ is the support of a numerically anticanonical divisor (cf. \cite{DlousskyNAC}).  In fact $T$ must be disconnected (\cite{AGG} Proposition 4), and by a result of Nakamura (\cite{FujikiPontASD} Lemma 3.3), a surface with disconnected numerical anticanonical divisor must be either a class $1$ Hopf surface, or a parabolic or hyperbolic Inoue surface.  The existence of generalized K\"ahler structures in some of these cases was established by Fujiki-Pontecorvo \cite{FujikiPontASD}.

\section{Generalized K\"ahler-Ricci flow} \label{s:GKRF}

As shown in (\cite{GKRF} Theorem 1.2), the pluriclosed flow preserves generalized K\"ahler geometry.  This arises as a consequence of Theorem \ref{t:PCFRGF}.  In particular, one notes that generalized K\"ahler structures consist of a metric which is pluriclosed with respect to two distinct complex structures, satisfying further integrability conditions.  By solving pluriclosed flow on both complex manifolds, and applying gauge transformations one obtains two solutions to (\ref{f:GRF}) with the same initial data, yielding a time dependent triple $(g_t, I_t, J_t)$ of generalized K\"ahler structures (cf. \cite{GKRF} for details).  Unpacking the construction yields the following evolution equations, which have the notable feature that the complex structures must evolve along the flow.

\begin{defn} A one-parameter family of generalized K\"ahler structures $(M^{2n}, g_t, I_t, J_t, H_t)$ is a solution of \emph{generalized K\"ahler-Ricci flow} (GKRF) if
\begin{gather}
\begin{split}
\dt g =&\ -2 \Rc^g + \frac{1}{2} H^2, \qquad \dt H = \gD_d H,\\
\dt I =&\ L_{\theta_I^{\sharp}} I, \qquad \dt J = L_{\theta_J^{\sharp}} J.
\end{split}
\end{gather}
\end{defn}

As a special case of pluriclosed flow, in principle we have already described all of the conjectural long time existence and convergence behavior in part II.  Nonetheless we will provide further discussion of long time existence results and refined descriptions of convergence behavior in this setting.  Interestingly, the GKRF in the $I$-fixed gauge preserves the Poisson structure on complex surfaces, a natural fact since it is holomorphic and thus rigid.  This can be shown using a case-by-case depending on the type (vanishing, nondegnerate, general) of the Poisson structure.  Conjecturally this is true in all dimensions, and should follow from a direct computation.  Depending on the type of the Poisson structure, the GKRF takes very different forms, which we will describe in turn below.

\subsection{Commuting case}

As discussed in \S \ref{s:GKG}, in the case $[I,J] = 0$ one obtains a splitting of the tangent bundle according to the eigenspaces of $Q = IJ$.  This also induces a splitting of the cotangent bundle, which induces a splitting of the exterior derivative
\begin{align*}
d = \del_+ + \delb_+ + \del_- + \delb_-.
\end{align*}
Arguing similarly to the $\del\delb$-Lemma in K\"ahler geometry, it is possible to obtain a local potential function describing generalized K\"ahler metrics in this setting.  In particular if $(g, I, J)$ is generalized K\"ahler then near any point there exists a smooth function $f$ such that
\begin{align*}
\gw_I = \i \left( \del_+ \delb_+ - \del_- \delb_- \right) f.
\end{align*}
This difference of sign indicates a fundamental distinction between K\"ahler geometry and generalized K\"ahler geometry in this setting: rather than being described locally by a plurisubharmonic function, the metric is described by a function which is plurisubharmonic in certain directions and plurisuperharmonic in others.  Nonetheless, given this local description, one expects the pluriclosed flow to reduce to a scalar PDE in this setting.  This was confirmed in (\cite{StreetsCGKFlow} Theorem 1.1), and locally this PDE takes the form
\begin{align} \label{f:TMA}
\frac{\del f}{\del t} = \log \frac{\det \i \del_+ \delb_+ f}{\det (- \i \del_- \delb_- f)}.
\end{align}
We refer to this equation as ``twisted Monge Ampere," as it is a natural combination of Monge Ampere operators for the different pieces of the metric.  This PDE is still parabolic, but is mixed concave/convex, and thus many of the usual methods for analyzing fully nonlinear PDE do not directly apply.  Nonetheless we were able to give a nearly complete picture of the long time existence of the flow on complex surfaces of this type.

\begin{thm} \label{t:CGKRF} (\cite{StreetsCGKFlow} Theorem 1.3) Let $(M^4, g_0, I, J)$ be a generalized K\"ahler
surface satisfying $[I,J] = 0$ and $I \neq \pm J$.  Suppose $(M^4,
I)$ is biholomorphic to
one of:
\begin{enumerate}
\item A ruled surface over a curve of genus $g \geq 1$.
\item A bi-elliptic surface,
\item An elliptic fibration of Kodaira dimension $1$,
\item A compact complex surface of general type, whose universal cover is
biholomorphic to $\mathbb H \times \mathbb H$,
\item An Inoue surface of type $S_M$.
\end{enumerate}
Then the solution to pluriclosed flow with initial condition $g_0$ exists on
$[0,\tau^*(\gw_0))$.
\end{thm}
Referring back to Theorem \ref{t:AG}, the only cases not covered by this theorem are ruled surfaces over curves of genus $0$ and Hopf surfaces.  The reason for the restriction in the theorem is that we require one of the line subbundles $T_{\pm}$ to have nonpositive first Chern class to obtain some partial a priori control over the metric.  Also we note that there is overlap between this theorem and Theorems \ref{t:negcurvthm}, \ref{t:torustheorem} above.

\subsection{Nondegenerate case} \label{ss:ndgk}

On the other extreme of generalized K\"ahler geometry is the nondegenerate case described above.  In \S \ref{s:GKG} we explained that if $(M^{4}, g, I, J, K)$ is hyperK\"ahler, we can interpret $(M^4, g, I, J)$ as a generalized K\"ahler structure.  Of course the underlying pairs $(g,I)$ and $(g,J)$ are still K\"ahler, but Joyce showed (cf. \cite{AGG}) that one can deform away from these examples to produce genuine examples of non-K\"ahler, generalized K\"ahler, structures with $\gs$ nondegenerate.  Specifically, given $f_t$ a one-parameter family of smooth functions, we define a family of vector fields $X_t$ via
\begin{align*}
X_t = \gs d f_t.
\end{align*}
Let $\phi_t$ denote the one parameter family of diffeomorphisms generated by $X_t$, which one notes is $\Omega$-Hamiltonian by construction (recall $\Omega = \gs^{-1}$).  The triple $(I, \phi_t^* J, \Omega)$ determines a generalized K\"ahler structure, with $g_t$ determined algebraically by (\ref{f:sigmadef}).

Surprisingly, the generalized K\"ahler-Ricci flow evolves by precisely this type of deformation.  In particular, given a generalized K\"ahler structure $(M^4, g, I, J)$ we let
\begin{align*}
p = \tfrac{1}{4} \tr IJ
\end{align*}
denote the angle between $I$ and $J$.  One has $\brs{p} \leq 1$, and the inequality is strict everywhere in the nondegenerate setting.  In four dimensions the function $p$ is constant if and only if $g$ is hyperK\"ahler (\cite{PontecorvoCS}).  Thus we expect the generalized K\"ahler-Ricci flow in this setting to completely reduce to quantities involving the angle function, and indeed this is the case.  Specifically, if we gauge-modify the generalized K\"ahler-Ricci flow to fix the complex structure $I$, then $J_t = \phi_t^* J$, where
\begin{align*}
\frac{d \phi}{d t} = \gs d \log \frac{1+p}{1-p}.
\end{align*}
This is roughly analogous to the usual scalar reduction for K\"ahler-Ricci flow, although here the Hamiltonian diffeomorphism $\phi_t$ depends on the entire history of the flow on $[0,t]$, and thus does not truly reduce to a single scalar.

We can give a complete description of the long time existence and a weak confirmation of the conjectured convergence behavior in this setting.
\begin{thm} (\cite{StreetsNDGKS} Theorem 1.1) Let $(M^4, g, I, J)$ be a nondegenerate generalized K\"ahler four-manifold.  The solution to generalized K\"ahler-Ricci flow with initial data $(g,I,J)$ exists for all time.  Moreover, $(\gw_I)_{t}$ subconverges to a closed current.
\end{thm}
\noindent While one expects smooth convergence to a hyperK\"ahler structure, the convergence behavior above at least shows that the flow contracts to the space of K\"ahler structures.  In fact more convergence properties can be shown (cf. \cite{StreetsNDGKS} for detail).
This result was extended to arbitrary dimensions by the author and Apostolov \cite{ASNDGKCY}.  The key observation is to show that a ``generalized Calabi-Yau quantity," motivated by natural constructions in generalized geometry (\cite{GualtieriThesis}), governs the dynamics of the flow in the same manner that $\log \frac{1+p}{1-p}$ does in four dimensions.

\subsection{General case}

The general case involves Poisson structures with type change locus, which exist on Del Pezzo surfaces, Hopf surfaces, and parabolic and hyperbolic Inoue surfaces as described above.  In all of these cases the conjectures on pluriclosed flow imply a connectedness result for the space of generalized K\"ahler structures.  Note that on a given complex manifold, the space of K\"ahler metrics is convex by linear paths, whereas uniqueness and moduli questions for complex structures are of course much more subtle.  In understanding the space of generalized K\"ahler metrics, these two problems are linked.  Moreover, there is no linear structure to this space, with the natural class of deformations instead using Hamiltonian diffeomorphisms as described above.
Thus the generalized K\"ahler-Ricci flow can potentially yield connectedness of the space of generalized K\"ahler structures, a nontrivial consequence due to the nonlinear structure of this space.  Furthermore, the conjectural behavior for parabolic and hyperbolic Inoue surfaces described in part II suggests that the limiting complete steady soliton metrics associated to these surfaces should in fact be generalized K\"ahler, and thus generalized K\"ahler structures should play a key role in understanding the geometrization of complex surfaces.

As an example, consider the case of $\mathbb {CP}^2$, where uniqueness of the complex structure in known \cite{YauCCann}.  Hitchin \cite{HitchindelPezzo} constructed generalized K\"ahler structures on $\mathbb{CP}^2$ by a modification of the Hamiltonian diffeomorphism method described above, deforming away from the standard Fubini Study structure.  Given this, and even knowing the uniqueness of complex structure on $\mathbb {CP}^2$, it is still possible that the space of generalized K\"ahler metrics has multiple disconnected components, i.e. there may exist other generalized K\"ahler triples not arising by this deformation.  This problem can be very naturally addressed using the generalized K\"ahler-Ricci flow, and in particular we can rule out the existence of such exotic generalized K\"ahler structures.  The main input is a description of the long time existence and weak convergence behavior of generalized K\"ahler-Ricci flow in this setting.

\begin{thm} \label{t:CP2} (\cite{SCPn}) Let $(\mathbb{CP}^2, g, I, J)$ be a generalized K\"ahler structure.  The solution to normalized generalized K\"ahler-Ricci flow with initial condition $(g,I,J)$ exists on $[0,\infty)$.  Moreover, $(\gw_I)_{t}$ subconverges to a closed current.
\end{thm}

As we detail in \cite{SCPn}, there is a natural completion of the space of generalized K\"ahler metrics extending the usual completion of K\"ahler metrics in the space of closed currents.  Our result yields connectivity of this space, and in fact that any point in this space is equivalent to the standard Fubini study structure by an extended Courant symmetry.  This is a natural extension of the classical uniqueness Theorem of Yau \cite{YauCCann} for complex structures on $\mathbb {CP}^2$ to generalized K\"ahler structures.

\end{document}